\let\accentvec\vec
\documentclass[smallextended,natbib]{svjour3}

\let\vec\accentvec
\usepackage{amsmath,amssymb,amsfonts,mathtools}
\usepackage{algorithm,algpseudocode}
\usepackage{graphicx}
\usepackage{multirow}
\usepackage{xspace}
\usepackage{enumitem}
\usepackage{rotating}
\usepackage{lmodern}
\usepackage[caption=false,subrefformat=parens,labelformat=parens]{subfig}
\usepackage{natbib}
\usepackage[naturalnames=false]{hyperref}
\algrenewcommand{\algorithmiccomment}[1]{\hfill$\blacktriangleright$ #1}

\hypersetup{colorlinks,
            linkcolor=blue,
            anchorcolor=blue,
            citecolor=blue}

\def\pdffigdir{figures-pdf}

\ifdefined\tikzdir
  \usepackage{pgfplots,tikzscale}
  \usetikzlibrary{external}
  \usetikzlibrary{arrows}
  \usetikzlibrary{calc}
  \usetikzlibrary{intersections}
  \pgfplotsset{compat=newest}
  \pgfplotsset{plot coordinates/math parser=false}
\else
  
\fi
\spnewtheorem{Asn}[theorem]{Assumption}{\bf}{\it}
\spnewtheorem{Thm}[theorem]{Theorem}{\bf}{\it}
\spnewtheorem{Lem}[theorem]{Lemma}{\bf}{\it}
\spnewtheorem{Cor}[theorem]{Corollary}{\bf}{\it}

\newcommand{\bseq}{\begin{subequations}}
\newcommand{\eseq}{\end{subequations}}

\newcommand{\be}{\begin{enumerate}}
\newcommand{\ee}{\end{enumerate}}
\newcommand{\bi}{\begin{itemize}}
\newcommand{\ei}{\end{itemize}}

\newcommand{\bR}{\mathbb{R}}    
\newcommand{\bOne}{\mathbf{1}}  

\newcommand{\cA}{{\cal A}}

\newcommand{\cD}{{\cal D}}

\newcommand{\cF}{{\cal F}}
\newcommand{\cG}{{\cal G}}

\newcommand{\cI}{{\cal I}}

\newcommand{\cL}{{\cal L}}

\newcommand{\cN}{{\cal N}}

\newcommand{\g}{\nabla}     

\DeclareMathOperator*{\s.t.}{subject\ to}
\DeclareMathOperator*{\tfor}{\ for\ }
\DeclareMathOperator*{\tif}{\ if\ }
\DeclareMathOperator*{\tow}{\ otherwise}
\DeclareMathOperator*{\dom}{dom}
\DeclareMathOperator*{\cone}{cone}
\DeclareMathOperator*{\cl}{cl}

\DeclareMathOperator*{\conv}{conv}

\newcommand{\lb}[1]{\underline{#1}}
\newcommand{\ub}[1]{\overline{#1}}

\newcommand{\omin}[1]{\min_{#1} \;\;}
\newcommand{\omax}[1]{\max_{#1} \;\;}

\newcommand{\ost}{\s.t. \;\;}

\DeclarePairedDelimiter\abs{\lvert}{\rvert}
\DeclarePairedDelimiter\norm{\lVert}{\rVert}
\DeclarePairedDelimiterX\innerp[2]{\langle}{\rangle}{#1,#2}
\providecommand\given{}
  \newcommand\SetSymbol[1][]{\nonscript\:#1\vert\nonscript\:
  \mathopen{}\allowbreak}
  \DeclarePairedDelimiterX\set[1]\{\}{%
  \renewcommand\given{\SetSymbol[\delimsize]}
  #1
}
\newcommand{\bm}[1]{\begin{bmatrix}#1\end{bmatrix}}

\def\<{\left<}
\def\>{\right>}
\def\({\left(}
\def\){\right)}
\def\mr{\multirow}
\def\mc{\multicolumn}

\def\noprint#1{}

\everymath{\displaystyle}

\newcommand{\tcellr}[2][c]{\begin{tabular}[#1]{@{}r@{}}#2\end{tabular}}
\newcommand{\tcell}[2][c]{\begin{tabular}[#1]{@{}c@{}}#2\end{tabular}}

\def\SLP{S$\ell_1$LP\xspace}
\def\SLP{S$\ell_1$LP\xspace}
\def\SLPAS{S$\ell_1$LP-AS\xspace}
\def\MATPOWER{MATPOWER\xspace}
\def\MATLAB{\textsc{Matlab}\xspace}
\def\IPOPT{\textsc{Ipopt}\xspace}
\def\CPLEX{\textsc{Cplex}\xspace}

\begin{document}

\journalname{Optimization and Engineering}

\title{An \SLP-Active Set Approach for Feasibility Restoration in
  Power Systems}

\author{Taedong~Kim \and Stephen~J.~Wright %
  \thanks{This research was supported by DOE Grant DE-SC000228, a DOE
  grant subcontracted through Argonne National Laboratory Award
  3F-30222, and National Science Foundation Grant DMS-1216318.}
}

\date{\today}

\institute{T. Kim \and S.J. Wright \at Computer Sciences Department, %
  1210 W. Dayton Street, University of Wisconsin,
  Madison, WI 53706, USA\\
  \email{tdkim@cs.wisc.edu, swright@cs.wisc.edu}}

\maketitle

\begin{abstract}
We consider power networks in which it is not possible to satisfy all
loads at the demand nodes, due to some attack or disturbance to the
network. We formulate a model, based on AC power flow equations, to
restore the network to feasibility by shedding load at demand nodes,
but doing so in a way that minimizes a weighted measure of the total
load shed, and affects as few demand nodes as possible.  Besides
suggesting an optimal response to a given attack, our approach can be
used to quantify disruption, thereby enabling ``stress testing'' to be
performed and vulnerabilities to be identified.
Optimization techniques
including nonsmooth penalty functions, sequential linear programming,
and active-set heuristics are used to solve this model. We describe an
algorithmic framework and present convergence results, including a
quadratic convergence result for the case in which the solution is
fully determined by its constraints, a situation that arises
frequently in the power systems application.
\end{abstract}

\keywords{AC power flow equations,
          composite nonsmooth optimization,
          sequential linear programming,
          active-set methods}

\section{Introduction}
\label{sec:intro}

Consider a power grid that has experienced an unexpected event that
may interfere with its ability to meet load requirements at its demand
nodes. The event could be a natural contingency, an equipment failure,
or a malicious attack intended to disrupt the system.  The operating
point needs to be modified after such an event. Since a power network
often has some resistance to small perturbations, the system may still
be operational after the disruption. If, however, the perturbation
exceeds the tolerance of the system, the power flow problem may fail
to have a solution until the configuration is adjusted, for example,
by reducing loads at the demand nodes.  We formulate and solve the
problem of shedding loads, in the least disruptive manner possible, to
restore feasible operation of the network. Our formulation is based on
the nonlinear AC power-flow model, so we need to solve a nonlinear
program, or minimize its equivalent nonsmooth penalty function.  The
measure of total load shed also serves as a metric for the severity of
the disruption, which may be useful in analyzing the vulnerability of
the grid to attacks of different kinds.

\subsection{Our Approach}
\label{sec:features}

Our approach for solving feasibility restoration problems based on AC
power-flow models uses an algorithmic framework of trust-region
methods for composite nonsmooth optimization (CNSO).  The subproblem
solved at each iteration can be posed as a linear program (LP), solved
with a simplex algorithm. Although the existing power-flow literature
tends to favor the use of interior-point methods, we find the LP-based
approach to be appealing because of its amenability to warm starting,
which often allows linearized subproblems to be solved quickly after
the first few ``outer'' iterations (of the sequential linear
programming strategy). Moreover, excellent software such as \CPLEX is
available for linear programming, and it can be invoked easily from
such power systems modeling frameworks such as
\MATPOWER~\citep{ZimMT11}. By contrast, warm-starting strategies for
interior-point methods have not proved to be effective in general
\citep{YilW02}, except when the optimal active set does not change
between outer iterations.

Because we formulate the nonlinear equality constraints in the problem
using an $\ell_1$-penalty term, rather than enforcing them as hard
constraints, all subproblems are feasible, provided that the initial
starting point for the very first subproblem is feasible.  Feasible
initial points can be chosen without any additional processing, since
the only explicit constraints left in the subproblems after
reformulation are box constraints. There is no need to complicate the
algorithm by solving separately for tangential and normal steps, as is
done in a number of previous approaches described below.

The first-order information used in sequential $\ell_1$-linear
programming (\SLP) framework is sometimes sufficient to produce a
quadratic convergence rate. When applied to a feasible system, our
approach usually takes the same steps as Newton's method applied to
the AC power flow equations (which are equality constraints in our
formulation), and both methods require only first derivatives of the
AC equations to be computed. Thus, if the power system is feasible,
the convergence rate of our algorithm is quadratic, like Newton's
method for nonlinear equations.  We observe this fast convergence
behavior too on several infeasible networks, especially when a network
can be recovered by adjusting just a few demand nodes.

In general, however, only linear convergence can be expected from
subproblems constructed with first-order information. To overcome this
slow rate of local convergence, we employ an active-set heuristic that
can accelerate the convergence of the algorithm at the cost of
second-derivative evaluations of the AC power flow formulae. In this
heuristic, we estimate the optimal active set of a problem and solve
the resulting equality-constrained optimization problem using a
Newton-like approach.  If the correct active set cannot be identified
after a modest number of attempts, we revert to the \SLP strategy.

\subsection{Previous Work}
\label{sec:prevwork}

There is a wide literature on applying nonlinear optimization
techniques to AC power flow equations. We survey here the works most
relevant to our approach from the recent literature, and then indicate
how our approach differs. Most of these papers use the setting of
optimal power flow (OPF) rather than feasibility, but since both these
problems are closely related to the formulation we consider here, we
discuss them together.

A formulation like ours for the unsolvability problem is considered in
\cite{GraMM96}, where the objective is a weighted sum of the fractions
of loads that cannot be met. A primal-dual interior-point method is
proposed for the resulting nonlinear program. The algorithm is a basic
interior-point method that includes few safeguards to ensure
convergence.

The formulation in \cite{BarS01a} seeks to minimize the sum of squares
of the loads that cannot be met, subject to the AC power flow
equations being satisfied at other specified nodes, along with line
and voltage magnitude limits. A primal-dual interior-point method is
applied to the resulting nonlinear program. In \cite{BarS01b}, the
same authors consider an alternative formulation in which the rate of
decrease of the load at the nodes eligible for load shedding is
specified in advance, leading to a nonlinear optimization formulation
whose objective consists of a single parameter. The interior-point
approach is again used to solve this formulation.

Trust-region approaches have been considered by several authors, in
conjunction with both interior-point and sequential quadratic
programming (SQP) algorithms.  All these papers have OPF as their
target problem, with the AC power flow equations as constraints, along
with line limits and voltage magnitude limits.  In \cite{ZhoZY05a}, an
SQP approach with trust regions and an $\ell_1$-merit function is
described. Since the primary trust-region subproblem can become
infeasible for small values of the trust-region radius $\Delta^k$, the
step is separated into components that are normal and tangent to the
active constraint manifold, respectively (see
\citet[Section~18.5]{NocW06} for a description of this technique), but
it is not clear how the active set is determined for purposes of this
calculation.  The second-order term in each subproblem is also not
specified.

\cite{SouT07} describe two trust-region approaches for the OPF
problem. The first approach is an SQP approach somewhat like that of
\cite{ZhoZY05a}, where the tangent subproblem is solved with an
interior-point method for quadratic programming. In the second
approach, a nonlinear primal-dual interior-point method is applied
directly to the nonlinear program. Second derivatives are used in both
approaches.  A journal paper by this team \citep{SouTC11a} focuses on
the first approach with some modifications, chiefly, that
interior-point methods for quadratic programming are used to solve
both the normal and tangential subproblems. Global convergence of the
approach to a local stationary point is noted, and computational
results are presented on standard test sets with up to 1211 buses.

The algorithm described in \cite{MinS05} is a sequential linear
programming (SLP) approach with trust regions, applied to an OPF
formulation, with the subproblems solved by an interior-point method
for linear programming. After the SLP step is calculated, the iterate
is adjusted by solving the AC power flow equations, thus ensuring that
every iterate satisfies the equality constraints exactly in the
nonlinear programming formulation. The authors note that the
trust-region subproblem can become infeasible if the radius $\Delta^k$
is too small, so they propose a modification in which the inequality
constraints are relaxed, and the amount of relaxation is penalized in
the subproblem. No convergence theory is presented.

\subsection{Outline}
\label{sec:outline}

In the next section, we introduce some notation along with the problem
formulation that is used as the basis of our analysis, together with
optimality conditions. Section~\ref{sec:slp} describes the sequential
$\ell_1$-linear programming (\SLP) framework and its global and local
convergence properties.  Section~\ref{sec:actset} discusses heuristics
for identifying the optimal active set and estimation of Lagrange
multipliers, which can be incorporated into the \SLP framework to
boost the convergence rate. Section~\ref{sec:formulation} describes
the application to minimal load shedding in disrupted power systems,
with computational results and comparisons to other approaches
described in Section~\ref{sec:result}.

\section{Formulation, Overview of the Algorithm, and Notation}
\label{sec:algoverview}

We begin this section by describing the problem formulation, along
with its optimality conditions and its nonsmooth penalty-function
equivalent. We also give a sketch of our algorithmic approach. The
second subsection presents basic relevant concepts from convex
analysis that are used in the convergence analysis of
Section~\ref{sec:slp}.

\subsection{Problem Formulations and Optimality Conditions}
\label{sec:formopt}

We consider the following optimization problem:
\bseq\label{eq:GNO} 
\begin{align}
  \omin{x} & p^Tx\\
  \ost & c(x) = 0                               \label{eq:GNO.eq}\\
       & \lb{x} \le x \le \ub{x},               \label{eq:GNO.ineq}
\end{align}
\eseq
where $p$, $\lb{x}$ and $\ub{x}$ are vectors in $\bR^n$ and
$c:\bR^n\rightarrow\bR^m$ is a nonlinear function with $n\ge m$.  Any
problem with nonlinear objectives and constraints can be formulated in
this way by introducing auxiliary variables; we prefer to isolate the
nonlinearity in the equality constraint to simplify the analysis and
description of the method.

At a feasible point $x$ of \eqref{eq:GNO}, we define the {\em set of
active inequalities} as
\bseq\label{eq:active}
\begin{align}
  \lb{\cA}(x) &= \set{i \given x_i = \lb{x}_i,\ i=1,2,\dotsc,n},\\
  \ub{\cA}(x) &= \set{i \given x_i = \ub{x}_i,\ i=1,2,\dotsc,n},\\
  \cA(x)      &= \lb{\cA}(x)\cup\ub{\cA}(x),
\end{align}
\eseq
and the {\em set of inactive inequalities} as
\bseq\label{eq:inactive}
\begin{align}
  \lb{\cI}(x) &= \set{1,2,\cdots,n}\backslash\lb{\cA}(x),\\
  \ub{\cI}(x) &= \set{1,2,\cdots,n}\backslash\ub{\cA}(x),\\
  \cI(x)      &= \set{1,2,\cdots,n}\backslash\cA(x).
\end{align}
\eseq
For convenience, we use abbreviated notation for these sets when
evaluated at an iterate $x^k$: $\lb{\cA}^k:=\lb{\cA}(x^k)$,
$\ub{\cA}^k:=\ub{\cA}(x^k)$, and so on. Similarly, at a solution $x^*$
of \eqref{eq:GNO}, we use $\lb{\cA}^* := \lb{\cA}(x^*)$, $\ub{\cA}^*
:= \ub{\cA}(x^*)$, and so on. The operator $| \cdot |$ is used to
denote cardinality of a set.

The {\em linear independent constraint qualification (LICQ)} for
\eqref{eq:GNO} at $x^*$ is:
\[
  \set{\g c_i(x^*),\ i=1,2,\dotsc,m} \cup \set{e_i,\, i\in\cA^*}
    \mbox{ is linearly independent}
\]
where $\g c_i(x^*)$ is the gradient of $c_i(x)$ at $x^*$ and $e_i$ is
the $i$th column of $n \times n$ identity matrix. We assume that LICQ
holds for the solutions of \eqref{eq:GNO}.

\begin{definition} \label{def:fully}
We say a solution $x^*$ of \eqref{eq:GNO} is {\it fully determined} by
the constraints if (a) $m+\abs{\cA^*}=n$; and (b) LICQ holds at
$x^*$. If $m+\abs{\cA^*}<n$, we call the solution {\it
underdetermined} (even if LICQ holds).
\end{definition}

The Lagrangian function $\cL(x,\lambda,\mu,\nu)$ of \eqref{eq:GNO} is
defined as
\[
  \cL(x,\lambda,\mu,\nu)
    := p^Tx+\lambda^Tc(x)-\mu^T(x-\lb{x})-\nu^T(x-\ub{x}),
\]
where $\lambda$, $\mu$, and $\nu$ are the Lagrange multipliers
corresponding to the equality constraints \eqref{eq:GNO.eq}, lower
bound and upper bound of \eqref{eq:GNO.ineq}, respectively.  The
first-order optimality conditions of \eqref{eq:GNO} are given as
\bseq \label{eq:kkt}
\begin{alignat}{1}
  \g_x\cL(x,\lambda,\mu,\nu) & = 0, \label{eq:kkt.L} \\
  c(x) & = 0,  \label{eq:kkt.c}  \\
  0 \le x - \lb{x} & \perp \mu \ge 0, \\
  0 \ge x - \ub{x} & \perp \nu \le 0,
\end{alignat}
\eseq
where $a\perp b$ indicates $a^Tb=0$ and
\begin{align*}
  \g_x\cL(x,\lambda,\mu,\nu) = p + \g c(x)^T\lambda-\mu-\nu.
\end{align*}

By introducing the following {\em $\ell_1$-penalty function} $\phi(x)$
with a penalty parameter $\omega>0$:
\begin{equation} \label{eq:CNSO}
  \phi(x) := p^Tx+\omega\norm{c(x)}_1,
\end{equation}
we can reformulate problem \eqref{eq:GNO} as a CNSO problem with  box
constraints:
\begin{equation} \label{eq:CNSO.constr}
  \omin{x} \phi (x) \;\; \ost \lb{x} \le x \le \ub{x}.
\end{equation}
It is well known that if the penalty parameter $\omega$ is
sufficiently large, under reasonable conditions, a local solution of
\eqref{eq:CNSO.constr} is a local minimizer of
\eqref{eq:GNO}~\cite[Theorem 17.3]{NocW06}. In our algorithm, the CNSO
problem is solved using a trust-region framework, for a particular
choice of parameter $\omega$. At each iteration, we define a
linearized model of the objective in \eqref{eq:CNSO.constr} and solve
an LP subproblem containing the constraints in \eqref{eq:CNSO.constr}
along with a trust region. Following \cite{Fle87}, we refer to this
approach as sequential $\ell_1$-linear programming (\SLP). The trust
region is adjusted so that the step obtained from the LP subproblem
gives a ``sufficient decrease'' in the objective $\phi$ at each
iteration, guaranteeing global convergence.  We will show that for
fully determined solutions, the \SLP algorithm converges quadratically
under certain conditions.

If the solution of a problem is {\em underdetermined}, we do not
expect the fast convergence of \SLP, so we enhance the basic strategy
with active-set heuristics that use second-order information, to
recover rapid local convergence.

Let $x^*$ be a (local) optimal solution of \eqref{eq:GNO}, satisfying
\eqref{eq:kkt}. Given the optimal active and inactive sets at $x^*$
(see \eqref{eq:active} and \eqref{eq:inactive}), we can rewrite the
conditions \eqref{eq:kkt} as follows:
\begin{alignat*}{2}
  \g_x\cL(x,\lambda,\mu,\nu) = 0, &\\
  c(x) = 0, &\\
  x_i-\lb{x}_i = 0, \quad \mu_i \ge 0 & \tfor i\in\lb{\cA}^*,\\
  x_i-\ub{x}_i = 0, \quad \nu_i \le 0 & \tfor i\in\ub{\cA}^*,\\
  \lb{x}_i < x_i < \ub{x}_i           & \tfor i\in\cI^*,\\
  \mu_i = 0                           & \tfor i\in\lb{\cI}^*,\\
  \nu_i = 0                           & \tfor i\in\ub{\cI}^*.
\end{alignat*}
By gathering the equality constraints in this system, we obtain
\begin{equation} \label{eq:KKT.eq}
  \bm{\g_x\cL(x,\lambda,\mu,\nu)\\
      c(x) \\
      (x-\lb{x})_{\lb{\cA}^*} \\
      (x-\ub{x})_{\ub{\cA}^*} \\
      \mu_{\lb{\cI}^*}\\
      \nu_{\ub{\cI}^*}}
  = 0,
\end{equation}
whose solution $(x^*,\lambda^*,\mu^*,\nu^*)$ also satisfies the
following inequalities:
\begin{equation} \label{eq:KKT.ineq}
\begin{aligned}
  \lb{x}_i < x_i < \ub{x}_i & \tfor i\in\cI^*,\\
  \mu_i \ge 0               & \tfor i\in\lb{\cA}^*,\\
  \nu_i \le 0               & \tfor i\in\ub{\cA}^*.
\end{aligned}
\end{equation}
If a starting point close enough to the optimum can be identified, and
if the optimal active sets are known, it may be possible to find a KKT
point \eqref{eq:kkt} by applying Newton's method for nonlinear
equations to \eqref{eq:KKT.eq}, then checking that the solution so
obtained satisfies \eqref{eq:KKT.ineq}.  To initiate this process, we
need reliable ways to identify the optimal active set, and to estimate
the values of the Lagrange multipliers. Under certain conditions, both
active sets and the optimal Lagrange multipliers can be estimated from
the duals of the LP subproblems at the previous \SLP iteration. If
after making a modest number of guesses of the optimal active set, the
strategy does not appear to be converging rapidly, we return to the
\SLP strategy.

\subsection{Convex Analysis Terminology and Notation}

Given a closed convex set $C \subset \bR^n$ we denote the {\em dual
cone} by $C^*$, where
\begin{equation} \label{eq:dualcone}
  C^* := \set{ p \given p^Tt \ge 0 \; \mbox{for all $t \in C$} }.
\end{equation}
The {\em polar cone} is denoted by $C^\circ$, where
\[
  C^\circ := \set{ p \given p^Tt \le 0 \; \mbox{for all $t \in C$} }.
\]
Note that $C^\circ = -C^*$.

The set $\Omega \subset \bR^n$ is a {\em polyhedral convex set} if
there is a finite collection of vectors $b_i \in \bR^n$ and scalars
$\gamma_i$, $i=1,2,\dotsc,K$, such that
\[
  \Omega = \set{ x \given b_i^T x \le \gamma_i, \; i=1,2,\dotsc,K }
\]
(see \cite[p.~170]{Roc70}). The active set $\mathcal{A}(x)$ at a given
$x \in \Omega$ is
\[
  \cA(x) := \set{ i = 1,2,\dotsc,K \given b_i^Tx = \gamma_i }.
\]
The normal cone $N_{\Omega}(x)$ to $\Omega$ at $x$ is defined as
\begin{equation} \label{eq:normalcone}
\begin{aligned}
N_{\Omega}(x)  := & \set{ y \given y^T(x'-x) \le 0, \; \mbox{for all $x'
  \in \Omega$} } \\
 = & \set*{ 0+\sum_{i \in \mathcal{A}(x)} \lambda_i b_i \given \lambda_i \ge 0 \; \mbox{for all $i \in \mathcal{A}(x)$} }.
\end{aligned}
\end{equation}
The set of feasible directions for $\Omega$ at $x \in \Omega$ is
defined as
\begin{equation} \label{eq:feasdir}
\begin{aligned}
\cF_{\Omega}(x) := & \set{ d \given x+ \alpha d\in\Omega
      \; \mbox{for all $\alpha$ sufficiently small and positive} } \\
= & \set{ d \given b_i^Td \le 0 \; \mbox{for all $i \in \mathcal{A}(x)$}}.
\end{aligned}
\end{equation}
It is easy to show (using a theorem of the alternative) that
$N_{\Omega}(x) = \cF_{\Omega}(x)^\circ$, which is the polar cone
constructed from $\cF_{\Omega}(x)$.

The subdifferential $\partial h(c)$ of a convex function $h: \bR^m \to
\bR$ is defined as
\[
\partial h(c) := \set{ v \given h(c) + v^T(c'-c) \le h(c') \;\;
                       \mbox{for all $c' \in \dom h$} }.
\]
When $h$ is a polyhedral convex function of the form
\begin{equation} \label{eq:polyc}
  h(c) := \omax{j=1,2,\dotsc,M} h_j^Tc + \beta_j,
\end{equation}
for some $h_j \in \bR^m$ and $\beta_j \in \bR$, $j=1,2,\dotsc,M$, we
obtain an explicit form of the subdifferential:
\[
\partial h(c) = \set*{ v = \sum_{j : h(c) = h_j^Tc + \beta_j} \lambda_j h_j \given
               \lambda_j \ge 0, \; \sum_{j : h(c) = h_j^Tc + \beta_j} \lambda_j = 1 },
\]
which in this case is a closed convex polyhedral
set~\cite[Theorem~19.1]{Roc70}.


\subsection{Other Notation}

We use $\bOne$ to denote the vector of ones, and $e_i$ to denote the
$i$th column of the identity matrix. Given a set $\cA \subset
\set{1,2,\dotsc,n}$, we use $I_{\cA}$ to denote the matrix whose rows
are the rows of the $n \times n$ identity matrix corresponding to the
entries in $\cA$.

\section{Sequential $\ell_1$-Linear Programming (S$\ell_1$LP)}
\label{sec:slp}

We begin this section by motivating and describing the \SLP approach
for problem \eqref{eq:CNSO.constr}.  Global convergence properties are
discussed in Subsection~\ref{sec:slp.global}, and local convergence
properties are the focus of Subsection~\ref{sec:slp.local}.

\subsection{Algorithm Description}

To solve \eqref{eq:CNSO.constr} using the \SLP framework, the LP
subproblems with an $\ell_{\infty}$-trust-region are defined as
follows, at iterate $x^k$:
\bseq \label{eq:SLP.subp}
\begin{align}
  \omin{d} & m^k(d) \\
  \ost     & \lb{x} \le x^k+d \le \ub{x} \label{eq:SLP.constr.x}\\
           & \norm{d}_\infty \le \Delta^k,
\end{align}
\eseq
where $\Delta^k > 0$ is a trust region radius and the linearized  {\em
model function} $m^k(d)$ is defined as follows:
\[
  m^k(d) := p^T(x^k+d) + \omega\norm{c(x^k)+\nabla c(x^k)d}_1.
\]
Note that $m^k(0)=\phi(x^k)$ by \eqref{eq:CNSO}. The function $m^k(d)$
is an approximation to the $\ell_1$-penalty function $\phi(x^k+d)$
that we ``trust'' to be a good approximation in the region
$\norm{d}_\infty \le \Delta^k$.  If a solution $d^k$ of
\eqref{eq:SLP.subp} yields a ``sufficient decrease'' in $\phi$, the
new iterate $x^{k+1}$ is define by $x^{k+1}=x^k+d^k$.  Otherwise, we
reject the step $d^k$, define the next iterate to be $x^{k+1}=x^k$,
reduce the trust-region radius, and proceed to the next iteration.

The quality of step $d^k$ is determined by means of the actual
reduction $\Delta\phi^k(d^k)$, the expected reduction $\Delta m^k(d^k)$,
and the ratio $\rho^k(d^k)$ between these quantities, defined as follows:
\bseq\begin{alignat}{2}
  \Delta\phi^k(d^k) &= \phi(x^k)-\phi(x^k+d^k),  & \quad
                  & \mbox{(Actual reduction)} \label{eq:SLP.red.act}\\
  \Delta m^k(d^k)   &= m^k(0)-m^k(d^k), &
                  & \mbox{(Expected reduction)} \label{eq:SLP.red.exp}\\
  \rho^k(d^k)       &= \frac{\Delta\phi^k(d^k)}{\Delta m^k(d^k)}. &
                  & \mbox{(Agreement ratio)} \label{eq:SLP.rho}
\end{alignat}\eseq
The expected reduction \eqref{eq:SLP.red.exp} for $d^k$ is always
nonnegative when $x^k$ is feasible for \eqref{eq:CNSO.constr}, since
$d=0$ is a feasible point of \eqref{eq:SLP.subp}. If $\rho^k(d^k)$
exceeds a positive threshold value $\lb{\rho}$, the step improves the
objective $\phi$ significantly, so we can accept it.  If $\rho^k(d^k)$
is close to 1, $m^k$ and $\phi$ are probably in good agreement over
the current trust region, so we increase the trust-region radius for
the next iteration. If $\rho^k(d^k)$ falls below another small
positive threshold $\lb{\eta}$, we deem the agreement between $m^k$
and $\phi$ to be poor, and we reduce the trust-region for the next
iteration. Since we choose parameters $\lb{\rho}$ and $\lb{\eta}$ to
satisfy $\lb{\rho}<\lb{\eta}$, the trust-region radius is always
reduced whenever the step is {\em not} taken.
The \SLP algorithm to solve \eqref{eq:CNSO.constr} is specified in
Algorithm~\ref{alg:SLP}.

\begin{algorithm}
\caption{Sequential $\ell_1$-Linear Programming}
\label{alg:SLP}
\begin{algorithmic}[1]
\Require
  \Statex Upper bound on trust-region radius $\ub{\Delta}>0$;
  \Statex Parameters $\lb{\rho}, \lb{\eta}, \ub{\eta}, c_1$, and $c_2$
    where $0<\lb{\rho}<\lb{\eta}<0.5<\ub{\eta}<1$, and $0<c_2<1<c_1$;
  \Statex Initial iterate $x^0$ and initial trust-region radius $\Delta^0\in(0,\ub{\Delta}]$;
\Ensure
  \Statex Solution $x^*$ of CSNO problem \eqref{eq:CNSO.constr};
\medskip
\For {$k=0,1,2,\cdots$}
  \State Construct and solve LP subproblem \eqref{eq:SLP.subp} at $x^k$ to obtain step $d^k$;
  \If{$\rho^k(d^k) < \lb{\rho}$}
    \Comment{Step $d^k$ is poor.}
    \State $x^{k+1}\gets x^k$;
  \Else
    \Comment{Step $d^k$ is good.}
    \State $x^{k+1}\gets x^k+d^k$;
  \EndIf
  \If {$\rho^k(d^k) < \lb{\eta}$}
    \Comment{Poor agreement.}
    \State $\Delta^{k+1}\gets c_2\Delta^k$;
  \ElsIf {$\rho^k(d^k) > \ub{\eta}$}
    \Comment{Good agreement.}
    \State $\Delta^{k+1}\gets \min\set{c_1\Delta^k, \ub{\Delta}}$;
  \Else
    \Comment{Fair agreement.}
    \State $\Delta^{k+1} \gets \Delta^k$;
  \EndIf
  \If {termination conditions are satisfied}
    \State {\bf break};
  \EndIf
\EndFor
\State $x^* \gets x^k$;
\end{algorithmic}
\end{algorithm}

Algorithm~\ref{alg:SLP} terminates if any of the following conditions
 are satisfied.
\begin{enumerate}[label=(\emph{\alph*})]
  \item\label{enum:term.1} The objective value is equal to the known
    best possible value of the problem.  (For example, in our
    application of Section~\ref{sec:formulation}, we can terminate if
    the load shedding is reduced to zero.)
  \item\label{enum:term.2} For some $\epsilon>0$, we have
    \begin{equation}
      \frac{\abs{\phi(x^k)-m^k(d^k)}}{\Delta^k}
        \le \epsilon.                       \label{eq:SLP.converged}
    \end{equation}
  \item\label{enum:term.3} The trust region radius drops below a lower
    bound $\lb{\Delta}$.
  \item\label{enum:term.4} The iteration counter reaches a maximum
    limit.
\end{enumerate}
Conditions \ref{enum:term.3} and \ref{enum:term.4} indicate
unsuccessful termination. If \ref{enum:term.2} is satisfied, the ratio
of possible improvement to trust region size is small, indicating that
the objective is in a flat region, probably near a solution.

In our implementation, we reformulate the
subproblem \eqref{eq:SLP.subp} as a true linear program by introducing
an auxiliary variables $\alpha$ to replace the $\ell_1$-norm, and
replacing the $\ell_{\infty}$-norm by bounds, as follows:
\bseq\label{eq:SLP.sublp}
\begin{align}
  \omin{d,\alpha} & p^Td+\omega\sum_{i=1}^m\alpha_i\\
  \ost & -\alpha \le c(x^k)+\nabla c(x^k)d \le \alpha\\
       & \max\set{\lb{x}-x^k,-\Delta^k}
         \le d \le
         \min\set{\ub{x}-x^k,\Delta^k},    \label{eq:SLP.LP.constr.d}
\end{align}
\eseq
where $\max$ and $\min$ in \eqref{eq:SLP.LP.constr.d} are applied
element-wise.  We have already mentioned that the subproblem is always
feasible no matter how small the trust-region radius $\Delta^k$ is,
provided $x^k$ is feasible for the simple constraints in
\eqref{eq:CNSO.constr}.  Due to the constraint
\eqref{eq:SLP.LP.constr.d}, the next iterate $x^{k+1}=x^k+d^k$ or
$x^{k+1}=x^k$ also satisfies the bounds in inequality constraints in
\eqref{eq:CNSO.constr}. Thus feasibility can be guaranteed for all
iterates $x^k$, provided that the initial iterate $x^0$ is feasible.
It follows that \eqref{eq:SLP.sublp} can be used in place of
\eqref{eq:SLP.subp} in Algorithm~\ref{alg:SLP} (line 2).

The \SLP approach has several advantages. First, it does not require
second derivative information, which may be expensive to evaluate.
Second, as we discuss later, rapid convergence can still be achieved
if the limit is a fully determined solution.  Third, the subproblem is
a linear program, so a simplex code can exploit warm start information
from the previous \SLP iteration.

\subsection{Global Convergence Properties}
\label{sec:slp.global}

In discussing the global convergence properties of \SLP we consider a
slight generalization of the formulation \eqref{eq:CNSO.constr}, which
we define as follows:
\begin{equation}\label{eq:GCNSO}
  \omin{x} \, \phi(x) := p^Tx+h(c(x)) \;\;
  \ost     x \in \Omega,
\end{equation}
where $\Omega\subseteq\bR^n$ is a polyhedral convex set and
$h:\bR^m\rightarrow\bR$ is a polyhedral convex function
\eqref{eq:polyc}.
The linear subproblem at iteration $k$ for \eqref{eq:GCNSO} is
\begin{equation}\label{eq:GCNSO.subp}
\begin{aligned}
  \omin{d} & m^k(d) := p^T(x^k+d) + h( c(x^k)+\nabla c(x^k)d )\\
  \ost     & x^k + d \in \Omega,\\
           & \norm{d}_\infty \le \Delta^k.
\end{aligned}
\end{equation}
We recover \eqref{eq:CNSO.constr} and \eqref{eq:SLP.subp},
respectively, by setting $h(x) = \omega \norm{x}_1$ and $\Omega =
\set{x \in \bR^n \given \lb{x} \le x \le \ub{x}}$ in \eqref{eq:GCNSO}
and \eqref{eq:GCNSO.subp}.

The following technical results are useful in proving global
convergence. (We omit the proof of the first result, which is
elementary.)
\begin{Lem}\label{lem:poly.dual}
  Let $C$ be a polyhedral convex set in $\bR^n$. Then the dual cone
  $C^*$ of $C$ (see \eqref{eq:dualcone}) is a closed convex polyhedral
  cone.
\end{Lem}

\begin{Lem}\label{lem:norm.vec}
  Let $C$ be a nonempty closed convex polyhedron and $D$ a nonempty
  closed convex polyhedral cone in $\bR^n$. If for each $d \in D$
  there exists $c \in C$ such that $c^Td \ge 0$, then there
  exists $c' \in C$ such that $(c')^T{d} \ge 0$ for all $d \in D$.
\end{Lem}
\begin{proof}
Note that $0 \in D$, since $D$ is a closed convex cone. If $0 \in C$,
we can choose $c'=0$ and we are done. If $0 \notin C$, assume for a
contradiction that for each $d \in D$ there exists $c\in C$ such that
$c^Td \ge 0$ but there is no $c' \in C$ such that $(c')^T{d} \ge 0$
for {\em all} $d \in D$. Let $D^*$ be the dual cone of $D$.  Then
$D^*$ is a closed convex polyhedral cone by Lemma \ref{lem:poly.dual},
and we have $C \cap D^* = \emptyset$. By the separating theorem for
polyhedra (see for example \citet[Theorem 10.4]{Van07}), there exists
a vector $p$ such that $p^Tc < 0$ for all $c \in C$ and $p^T{d^*} \ge
0$ for all $d^* \in D^*$. From the latter, we have $p \in D^{**} =
\cl\conv\cone D = D$. Thus we have identified $p \in D$ such that
$p^Tc < 0$ for all $c \in C$. This contradicts the assumption
mentioned before. Thus there exists $c' \in C$ such that $(c')^Td \ge
0$ for all $d \in D$. \qed
\end{proof}

Now we are ready to show the global convergence of the \SLP algorithm.
For this purpose we recall the definition \eqref{eq:normalcone} of the
normal cone $N_{\Omega}(x)$ to $\Omega$ at a point $x \in \Omega$.
We also need the first-order optimality condition for \eqref{eq:GCNSO}
at a point $x^*$, which is that there exists $\lambda^* \in \partial
h(c(x^*))$ such that
\begin{equation} \label{eq:GCNSO.kkt}
0 \in p + \nabla c(x^*)^T\lambda^* + N_\Omega(x^*).
\end{equation}

\begin{Thm}\label{thm:glob.conv}
  Assume that $\phi(x)$ is bounded below and the sequence $\set{x^k}$
  generated by Algorithm~\ref{alg:SLP} is bounded. Also assume that
  $\partial h (c)$ is bounded for all $c$.  Then there is a
  subsequence $S$ of $\set{x^k}$ with a limit point $x^\infty$ that
  satisfies the first-order optimality conditions
  \eqref{eq:GCNSO.kkt}.
\end{Thm}
We note that the assumptions for Theorem~\ref{thm:glob.conv} hold for
our problem since $h(c)=\omega\norm{c}_1$ and so $\partial h(x)
\subset \otimes_{i=1}^m [-\omega,\omega]$, while the objective
$\phi(x)$ is bounded below by zero.  Versions of this result without
the constraint $x \in \Omega$ appear in \citet[Theorem 14.5.1]{Fle87}
and \citet[Theorem 2.1]{FleS89}. Our proof follows these, with
modifications to handle the presence of the feasible set $\Omega$.
\begin{proof}
By taking a further subsequence $S$ if necessary, we can assume that
one of the following two cases occurs for indices $k \in S$:
\begin{enumerate}[label=(\emph{\alph*})]
    \item\label{enum:case1} $\rho^k(d^k) < \lb{\eta}$ and $\Delta^{k+1}
      \rightarrow 0$.  Thus $\norm{d^k}\rightarrow0$.
    \item\label{enum:case2} $\rho^k(d^k) > \lb{\eta}$ and $\inf \Delta^k > 0$.
\end{enumerate}
We consider first case \ref{enum:case1}.  Suppose for contradiction
that there exists a nonzero vector $s\in \cF_\Omega(x^\infty)$ and a
constant $\beta>0$ such that
\begin{equation}\label{eq:beta}
 \omax{ \lambda\in\partial h( c(x^\infty) ) }
  \frac{1}{\norm{s}} s^T (p+\nabla c(x^\infty)^T\lambda)
    = -\beta, \quad \beta>0.
\end{equation}
where the closed convex cone of feasible directions $\cF_\Omega(x)$ is
defined in \eqref{eq:feasdir}. Because $\Omega$ is polyhedral convex,
we have $N_\Omega(x) \subseteq N_\Omega(x^\infty)$ for all $x$ close
enough to $x^\infty$, by outer-semicontinuity of $N_\Omega(x)$
\citep[Proposition 6.6]{RocW09}. (This fact is a consequence of
$\cA(x) \subset \cA(x^\infty)$ for $x$ close enough to $x^\infty$.)
We thus have $\cF_\Omega(x^\infty)\subseteq\cF_\Omega(x)$, since
$\cF_\Omega(x) = N_\Omega(x)^\circ$. Since $s \in
\cF_\Omega(x^\infty)$ and $\lim_{k \in S} x^k=x^\infty$, we have $s
\in \cF_\Omega(x^k)$ for $k$ sufficiently large.
Furthermore, it is easy to show from the definitions associated with
polyhedral convex $\Omega$ at the end of Section~\ref{sec:algoverview}
that $x^k + \norm{d^k}({s}/{\norm{s}}) \in \Omega$, since $\norm{d^k}
\rightarrow 0$ and $x^k \rightarrow x^\infty$. Since $c(x)$ is
continuous and $\partial h$ is bounded, we know from Taylor's Theorem
that
\begin{equation}\label{eq:linear.error}
  \begin{aligned}
    \phi(x^k+d) &= p^T(x^k+d) + h(c(x^k+d)) \\
          &= p^T(x^k+d) + h\(c(x^k) + \nabla c(x^k)d+o(\norm{d})\)\\
          &= p^T(x^k+d) + h(c(x^k) + \nabla c(x^k)d)+o(\norm{d})\\
          &= m^k(d)+o(\norm{d}).
  \end{aligned}
\end{equation}
Then for sufficiently large $k$, we have
\begin{align*}
  \Delta m^k(d^k)
    &=   \phi(x^k)-m^k(d^k) \\
    &\ge \phi(x^k)-m^k\(\norm{d^k}\frac{s}{\norm{s}}\)\\
    &=   \phi(x^k)-\phi\(x^k+\norm{d^k}\frac{s}{\norm{s}}\)+o\(\norm{d^k}\)\\
    &\ge \beta\norm{d^k}+o\(\norm{d^k}\).
\end{align*}
The first inequality follows from the fact that $d^k$ solves the
subproblem \eqref{eq:GCNSO.subp} whereas $s (\norm{d^k}/\norm{s})$ is
another feasible point for this subproblem. The second equality
follows from \eqref{eq:linear.error}, and the final inequality follows
from \citet[Corollary to Lemma 14.5.1]{Fle87} and \eqref{eq:beta}.
Also we know that
\[
  \Delta\phi^k(d^k) = \Delta m^k(d^k)+o(\norm{d^k}),
\]
so that
\[
  \rho^k(d^k) = \frac{\Delta\phi^k(d^k)}{\Delta m^k(d^k)} = 1 + o(1).
\]
This contradicts the assumption that $\rho^k(d^k) < \lb{\eta}$.
Therefore there is no such vector $s\in\cF_\Omega(x^\infty)$ for which
\eqref{eq:beta} holds, so that for all $s\in\cF_\Omega(x^\infty)$, we
have
\[
  \omax{\lambda\in\partial h(c(x^\infty))}
  \frac{1}{\norm{s}} s^T (p+\nabla c(x^\infty)^T\lambda)  \ge 0.
\]
We can now set $D = \cF_\Omega(x^\infty)$ and $C = \partial
h(c(x^{\infty}))$ in Lemma \ref{lem:norm.vec} to conclude that there
exists $\lambda^* \in \partial h(c(x^\infty))$ such that
$s^T{(p+\nabla c(x^\infty)^T\lambda^*)} \ge 0$ for all $s \in
\cF_\Omega(x^\infty)$. We therefore have the desired result that $0
\in p + \nabla c(x^\infty)^T\lambda^* + N_\Omega(x^\infty)$.

In case \ref{enum:case2}, we know that $\Delta\phi^k(d^k) \rightarrow
0$ from
\[
  \phi(x^0) - \phi(x^\infty) \ge \sum_{k\in S} \Delta\phi^k(d^k).
\]
Since $\rho^k(d^k) \ge \lb{\eta}$, we also have $\Delta m^k(d^k)
\rightarrow 0$.  Now let $d^\infty$ be a solution of the subproblem
\eqref{eq:GCNSO.subp} with $x^k=x^\infty$ and $\Delta^\infty <
\tilde{\Delta} := \inf_{k\in S}\Delta^k$. Also define $\tilde{x} =
x^\infty+d^\infty$. Then
\[
  \norm{\tilde{x}-x^k}
     \le \norm{\tilde{x}-x^\infty} + \norm{x^\infty-x^k}
     =   \norm{d^\infty}+o(1)
     \le \Delta^\infty+o(1)
     \le \Delta^k,
\]
for sufficiently large $k \in S$. Thus $\tilde{x}-x^k$ is feasible for
\eqref{eq:SLP.subp}
and
\[
  m^k(\tilde{x}-x^k) \ge m^k(d^k) = \phi(x^k) - \Delta m^k(d^k).
\]
By taking limits of both sides, we have
\[
  m^\infty(d^\infty) \ge \phi(x^\infty) = m^\infty(0).
\]
Therefore, $0$ is also a solution of \eqref{eq:GCNSO.subp} (with
$x^k =x^\infty$ and $\Delta^k = \Delta^\infty$), in particular, the
trust-region constraint is inactive. From the optimality conditions
of \eqref{eq:GCNSO.subp}, there exists $\lambda^* \in \partial
h(c(x^\infty))$ such that
\[
  0 \in p + \nabla c(x^\infty)^T\lambda^* + N_\Omega(x^\infty).
\]
Thus $x^\infty$ is a KKT point of \eqref{eq:GCNSO}.  \qed
\end{proof}

\subsection{Fast Local Convergence in the Fully Determined Case}
\label{sec:slp.local}

As discussed in Section~\ref{sec:algoverview}, we cannot expect the
\SLP algorithm to have a local convergence rate faster than linear in
general, because it uses only first-order information about the
function $c(x)$. However, when the solution $x^*$ is fully determined
(see Definition~\ref{def:fully}), and when certain other conditions
hold, the algorithm converges locally at a quadratic rate. In this
subsection, we analyze this phenomenon, which is commonly observed in
our application of Section~\ref{sec:formulation}.

First-order optimality conditions for problem \eqref{eq:CNSO.constr}
are as follows (cf. \eqref{eq:kkt}):
\bseq\label{eq:CNSO.KKT}
\begin{align}
  p + \omega\nabla c(x)^T\lambda-\mu-\nu &= 0, \label{eq:CNSO.KKT.lagr}\\
  \lambda & \in \partial\norm{c(x)}_1, \\
  0 \le x-\lb{x} & \perp \mu \ge 0, \\
  0 \ge x-\ub{x} & \perp \nu \le 0,
\end{align}
\eseq
where $\lambda, \mu$ and $\nu$ are the dual variables. We say that
{\em strict complementarity} holds at a primal-dual solution
$(x^*,\lambda^*,\mu^*,\nu^*)$ of \eqref{eq:CNSO.KKT} if the following
conditions are satisfied:
\begin{enumerate}[label=(\emph{\alph*})]
  \item $\lambda_j^* \in (-1,1)$ if $c_j(x^*)=0$ for $j \in\set{1,2,\dotsc,m}$,
  \item $\mu_i^* > 0 \tif (x^*-\lb{x})_i = 0 \tfor
    i\in\set{1,2,\dotsc,n}$, that is, for $i \in \lb{\cA}^*$,
  \item $\nu_i^* < 0 \tif (x^*-\ub{x})_i = 0 \tfor
    i\in\set{1,2,\dotsc,n}$, that is, for $i \in \ub{\cA}^*$.
\end{enumerate}
In other words, strict complementarity requires existence of
$\gamma>0$ such that
\bseq \label{eq:sc.gamma}
\begin{alignat}{2}
  \abs{\lambda_j^*} & \le 1-\gamma \;\; && \mbox{for $j \in\set{1,2,\dotsc,m}$}, \\
  \mu^*_i &\ge \gamma \;\; && \mbox{for $i \in \lb{\cA}^*$}, \\
  \nu^*_i &\le -\gamma \;\; && \mbox{for $i \in \ub{\cA}^*$}.
\end{alignat}
\eseq

Guided by the optimality conditions \eqref{eq:CNSO.KKT}, we have the
following result concerning solution $d^k$ of the LP subproblem
\eqref{eq:SLP.sublp}.
\begin{Lem}\label{lem:local.conv}
  Suppose that $(x^*,\lambda^*,\mu^*,\nu^*)$
  satisfies \eqref{eq:CNSO.KKT}, that $x^*$ is fully determined, and
  that strict complementarity holds at
  $(x^*,\lambda^*,\mu^*,\nu^*)$. Then for all $x^k$ sufficiently close
  to $x^*$, and provided that $\Delta^k \ge 2 \| x^k-x^* \|$, the
  subproblem \eqref{eq:SLP.subp} has a solution $d^k$ such that
  $c(x^k)+\nabla c(x^k)d^k = 0$, $(x^k+d^k)_i=\lb{x}_i$ for
  $i \in \lb{\cA}^*$, $(x^k+d^k)_i=\ub{x}_i$ for $i \in \ub{\cA}^*$
  and $\norm{x^k+d^k-x^*}=O(\norm{x^k-x^*}^2)$.
\end{Lem}
\begin{proof}
We write the first-order necessary conditions of the subproblem
\eqref{eq:SLP.subp} {\em without} the trust-region constraints as
follows:
\bseq\label{eq:SLP.KKT}
\begin{align}
  p + \omega\nabla c(x^k)^T\lambda-\mu-\nu &= 0     \label{eq:SLP.KKT.lagr}\\
  \lambda &\in \partial\norm{c(x^k)+\nabla c(x^k)d}_1 \\
  0 \le x^k+d-\lb{x} &\perp \mu \ge 0\\
  0 \ge x^k+d-\ub{x} &\perp \nu \le 0.
\end{align}
\eseq
We will construct a solution $(d^k,\lambda^k,\mu^k,\nu^k)$ to these
conditions and then show that $d^k = x^*-x^k+O(\norm{x^k-x^*}^2)$ and
thus $\norm{d^k} \le 2 \|x^k-x^* \| \le \Delta^k$ when
$\norm{x^k-x^*}$ is sufficiently small. In addition, this $d^k$
satisfies the other conditions mentioned in the theorem, and together
with $c(x^k)+\nabla c(x^k)d^k=0$, it solves \eqref{eq:SLP.subp}.

In constructing our solution $(d^k,\lambda^k,\mu^k,\nu^k)$ to
\eqref{eq:SLP.KKT}, we define
\bseq \label{eq:SLP.munu}
  \begin{align}
    \mu_i^k = 0 &\tfor i \notin \lb{\cA}^*\\
    \nu_i^k = 0 &\tfor i \notin \ub{\cA}^*.
  \end{align}
\eseq
and
\bseq\label{eq:SLP.constr.ineq}
  \begin{align}
    d_i^k = \lb{x}_i-x_i^k = x_i^*-x_i^k &\tfor i \in \lb{\cA}^*\\
    d_i^k = \ub{x}_i-x_i^k = x_i^*-x_i^k &\tfor i \in \ub{\cA}^*.
  \end{align}
\eseq
We further require $d^k$ to satisfy
\begin{equation} \label{eq:cd}
    c(x^*) = 0 = c(x^k)+\nabla c(x^k)d^k.
\end{equation}
The remaining components of  $(d^k,\lambda^k,\mu^k,\nu^k)$ are
required to satisfy the following linear system:
\begin{equation} \label{eq:linsys}
    \bm{\nabla c(x^k)^T &
        I_{\lb{\cA}^*}^T &
        I_{\ub{\cA}^*}^T}
    \bm{ \omega\lambda^k\\
         -\mu_{\lb{\cA}^*}^k\\
         -\nu_{\ub{\cA}^*}^k}
    = -p,
\end{equation}
which is essentially \eqref{eq:SLP.KKT.lagr} with the substitution
\eqref{eq:SLP.munu}. By comparing \eqref{eq:linsys} with
\eqref{eq:CNSO.KKT.lagr}, we have
\begin{align*}
  \bm{\nabla c(x^k)^T &
      I_{\lb{\cA}^*}^T &
      I_{\ub{\cA}^*}^T}
  \bm{ \omega\lambda^k\\
       -\mu_{\lb{\cA}^*}^k\\
       -\nu_{\ub{\cA}^*}^k}
  &  =
      \bm{\nabla c(x^*)^T &
          I_{\lb{\cA}^*}^T &
          I_{\ub{\cA}^*}^T}
      \bm{ \omega\lambda^*\\
          -\mu_{\lb{\cA}^*}^*\\
          -\nu_{\ub{\cA}^*}^*}\\
  &  =
     \bm{\nabla c(x^k)^T &
         I_{\lb{\cA}^*}^T & I_{\ub{\cA}^*}^T}
     \bm{ \omega\lambda^*\\
          -\mu_{\lb{\cA}^*}^*\\
          -\nu_{\ub{\cA}^*}^*}+O(\norm{x^k-x^*}).
\end{align*}
Since LICQ holds at $x^*$ and $m+\abs{\cA^*}=n$, $\bm{\nabla c(x^*)\\
I_{\cA^*}}$ is a nonsingular square matrix. Then by continuity of
$\nabla c$, the matrix $\bm{\nabla c(x^k)\\I_{\cA^*}}$ is uniformly
nonsingular for $\norm{x^k - x^*}$ sufficiently small, and we have
\[
  \bm{ \omega\lambda^k\\
      -\mu_{\lb{\cA}^*}^k\\
      -\nu_{\ub{\cA}^*}^k}
  = \bm{ \omega\lambda^*\\
        -\mu_{\lb{\cA}^*}^*\\
        -\nu_{\ub{\cA}^*}^*}
    + O(\norm{x^k-x^*}).
\]
It follows from this relation and the strict complementarity
assumption on $(x^*,\lambda^*,\mu^*,\nu^*)$ that $\mu_{\lb{\cA}^*}^k >
0$, $\nu_{\ub{\cA}^*}^k < 0$, and $\lambda^k \in (-1,1)^m$ for
$\norm{x^k-x^*}$ sufficiently small. Because of \eqref{eq:cd}, these
values are feasible for \eqref{eq:SLP.KKT}.

Note from  \eqref{eq:cd}  and Taylor's theorem that
\[
  \nabla c(x^k)d^k
  = c(x^*)-c(x^k)
  = \nabla c(x^k)(x^*-x^k) + O(\norm{x^*-x^k}^2).
\]
By combining this expression with \eqref{eq:SLP.constr.ineq}, we
obtain
\[
  \bm{ \nabla c(x^k)\\ I_{\cA^*} }d^k
  = \bm{ \nabla c(x^k)\\ I_{\cA^*} }(x^*-x^k)
    + O(\norm{x^*-x^k}^2).
\]
By invertibility of the coefficient matrix (discussed above), we have
\begin{equation} \label{eq:SLP.d}
  d^k = (x^*-x^k) + O(\norm{x^*-x^k}^2).
\end{equation}
This estimate implies that the bounds inactive at $x^*$ will also be
inactive at $x^k+d^k$, where $d^k$ is the solution of the subproblem
\eqref{eq:SLP.subp}, for $x^k$ close enough to $x^*$.

At this point, we have found $(d^k,\lambda^k,\mu^k,\nu^k)$ that
satisfies \eqref{eq:SLP.KKT}, and that $d^k$ satisfies the other
properties claimed in the theorem.  Moreover, $d^k$ satisfies the
trust-region bound, since
\[
\| d^k \| = \| x^k-x^*\| + O(\|x^k-x^*\|^2) \le 2 \|x^k-x^* \| \le \Delta^k,
\]
for $x^k$ close enough to $x^*$. Thus $d=d^k$ solves
\eqref{eq:SLP.subp}, as required.  \qed
\end{proof}

Lemma~\ref{lem:local.conv} says that under the given conditions, there
exists a solution $d^k$ of the linearized subproblem that is a
quadratic step to a solution of the problem \eqref{eq:CNSO.constr}. We
show now that $d^k$ is {\em unique} solution of the linearized
subproblem, under the given conditions.
\begin{Lem}
  Under the conditions of Lemma \ref{lem:local.conv}, the vector $d^k$
  described in this result is the unique
  solution of the linearized subproblem \eqref{eq:SLP.sublp}.
\end{Lem}
\begin{proof}
Assume that the conditions in Lemma \ref{lem:local.conv} hold and let
$(d^k, \lambda^k, \mu^k, \nu^k)$ satisfy \eqref{eq:SLP.KKT}.  From
\eqref{eq:SLP.KKT.lagr}, we have for all $d$ that
\[
  (p+\omega\g c(x^k)^T\lambda^k)^T(d-d^k) = (\mu^k+\nu^k)^T(d-d^k).
\]
Since $m^k(d)$ is a convex function, we have
\begin{align*}
  m^k(d) &\ge m^k(d^k) + (p+\omega\g c(x^k)^T\lambda^k)^T(d-d^k) \\
         &=   m^k(d^k) + (\mu^k+\nu^k)^T(d-d^k).
\end{align*}
We also know that
\begin{equation}\label{eq:SLP.KKT.lagr.sign}
  (\mu^k+\nu^k)_i
  \begin{cases}
    > 0 & \tif i \in \lb{\cA}^*, \\
    < 0 & \tif i \in \ub{\cA}^*, \\
    = 0 & \tow.
  \end{cases}
\end{equation}

Now assume for contradiction that $d'$ is a solution
of \eqref{eq:SLP.subp}, different from $d^k$. If $d'_{\cA^*} =
d^k_{\cA^*}$, then both of $d'$ and $d^k$ satisfy the following
equation:
\[
  \bm{\nabla c(x^k)\\
      I_{\cA^*}}
  d
  = \bm{c(x^*)-c(x^k)\\
        (x^*-x^k)_{\cA^*}},
\]
which implies that $d'=d^k$, by nonsingularity of the coefficient
matrix, contradicting the choice of $d'$.  If $d'_{\cA^*} \neq
d^k_{\cA^*}$, we have since $d^k_i$ achieves its lower bound for
$i \in \lb{\cA}^*$ and its upper bound for $i \in
\ub{\cA}^*$ that
\[
    (d'-d^k)_{\lb{\cA}^*} \ge 0, \quad
    (d'-d^k)_{\ub{\cA}^*} \le 0, \quad
    (d'-d^k)_{\cA^*} \neq 0.
\]
By \eqref{eq:SLP.KKT.lagr.sign}, we therefore have that
$(\mu^k+\nu^k)^T(d'-d^k)>0$, which implies
\[
    m^k(d')
    \ge m^k(d^k)+(\mu^k+\nu^k)^T(d'-d^k)
   > m^k(d^k).
\]
Thus $d'$ cannot be a solution of \eqref{eq:SLP.subp}. \qed
\end{proof}

We now show quadratic local convergence of Algorithm~\ref{alg:SLP}
applied to \eqref{eq:CNSO.constr}, under the assumptions of this
section.
\begin{Thm}
Suppose that Algorithm~\ref{alg:SLP} generates a sequence $\{x^k\}$
  that has an accumulation point $x^*$, where $x^*$
  satisfies \eqref{eq:CNSO.KKT}, is fully determined, and 
  strict complementarity holds at $(x^*,\lambda^*,\mu^*,\nu^*)$.
  Suppose for some iterate $\bar{k}$ with $\|x^{\bar{k}}-x^*\|$
  sufficiently small, we have that $\Delta^{\bar{k}} \ge 2 \|
  x^{\bar{k}}-x^* \|$. Then the sequence $\{ x^k \}$ converges
  quadratically to $x^*$.
\end{Thm}
\begin{proof}
The main part of the proof is to show that the step $d^k$ defined in
Lemma~\ref{lem:local.conv} is accepted by Algorithm~\ref{alg:SLP},
with a value of $\rho^k(d^k)$ close to $1$, for $\|x^k-x^*\|$
sufficiently small with $\Delta^k \ge 2 \|x^k-x^* \|$. It follows that
the trust region radius is not decreased, and hence remains inactive
at the next iteration. A recursive argument applied to the iteration
sequence starting at $k=\bar{k}$ completes the proof.

For $d^k$ defined in Lemma~\ref{lem:local.conv}, we have the definition
\eqref{eq:SLP.rho} that
\begin{align*}
  \rho^k(d^k)
         &= \frac{\phi(x^k)-\phi(x^k+d^k)}{m^k(0)-m^k(d^k)} \\
         &= \frac{-p^Td^k+\omega(\norm{c(x^k)}_1-\norm{c(x^k+d^k)}_1)}
           {-p^Td^k+\omega(\norm{c(x^k)}_1-\norm{c(x^k)+\nabla c(x^k)d^k}_1)} \\
         &= 1+\frac{-\omega\norm{c(x^k+d^k)}_1}
                   {-p^Td^k+\omega\norm{c(x^k)}_1} \quad\quad
  \mbox{since $c(x^k)+\nabla c(x^k)d^k=0,$} \\
         &= 1+\frac{-\omega\norm{c(x^k)+\nabla c(x^k)d^k+O(\norm{d^k}^2)}_1}
                   {-p^Td^k+\omega\norm{\nabla c(x^k)d^k}_1}\\
         &= 1+\frac{O(\norm{d^k}^2)}{-p^Td^k+\omega\norm{\nabla c(x^k)d^k}_1}.
\end{align*}
We show now that the remainder term in this expression is of size
$O(\norm{d^k})$ by showing that the denominator is bounded below by a
multiple of $\norm{d^k}$.  From the optimality condition
\eqref{eq:CNSO.KKT.lagr}, we have
\begin{equation}
  p^Td^k + \omega (\lambda^*)^T\g c(x^*)d^k
  = (\mu^*)^Td^k+(\nu^*)^Td^k.            \label{eq:CNSO.KKT.lagr.eq}
\end{equation}
Since strict complementarity holds at $(x^*,\lambda^*,\mu^*,\nu^*)$,
we have from the constant $\gamma>0$ defined in \eqref{eq:sc.gamma}
that
\[
  (\lambda^*)^T \nabla c(x^*) d^k
    \ge -(1-\gamma) \norm{ \nabla c(x^*)d^k }_1
\]
and so
\begin{equation}
  \norm{\g c(x^*)d^k}_1
  \ge \gamma \norm{\g c(x^*)d^k}_1
      - (\lambda^*)^T\g c(x^*)d^k.             \label{eq:lambda.ineq}
\end{equation}
Using \eqref{eq:sc.gamma} again, noting
that $\mu^*_i \ge \gamma$ and $d_i^k \le 0$ for $i \in \lb{\cA}^*$,
and $\nu_i^* \le -\gamma$ and $d_i^k \ge 0$ for $i \in \ub{\cA}^*$, we
also have
\begin{equation} \label{eq:zeta.new}
  -(\mu_{\lb{\cA}^*}^*)^Td_{\lb{\cA}^*}^k
  -(\nu_{\ub{\cA}^*}^*)^Td_{\ub{\cA}^*}^k
  \ge \gamma \norm{d^k_{\cA^*}}_1.
\end{equation}
By combining \eqref{eq:SLP.d}, \eqref{eq:CNSO.KKT.lagr.eq}, \eqref{eq:lambda.ineq}, and
\eqref{eq:zeta.new}, we have
\begin{align*}
 -p^Td^k & +\omega\norm{\g c(x^k)d^k}_1 \\
      &=   -p^Td^k+\omega\norm{\g c(x^*)d^k}_1 + O(\norm{d^k}^2) \\
      &\ge -p^Td^k
           -\omega(\lambda^*)^T\g c(x^*)d^k+\omega\gamma\norm{\g c(x^*)d^k}_1
           + O(\norm{d^k}^2)\\
      &=   -(\mu^*)^Td^k-(\nu^*)^Td^k
           +\omega\gamma\norm{\g c(x^*)d^k}_1
           + O(\norm{d^k}^2) \\
      &=   -(\mu_{\lb{\cA}^*}^*)^Td_{\lb{\cA}^*}^k
           -(\nu_{\ub{\cA}^*}^*)^Td_{\ub{\cA}^*}^k
           +\omega\gamma\norm{\g c(x^*)d^k}_1
           + O(\norm{d^k}^2)\\
      &\ge  \gamma\norm{d^k_{\cA^*}}_1
           +\omega\gamma\norm{\g c(x^*)d^k}_1
           + O(\norm{d^k}^2)\\
      &=    \gamma\norm*{\bm{\omega\g c(x^*)\\
                             I_{\cA^*}}d^k}_1
           + O(\norm{d^k}^2)\\
      &\ge  \gamma\zeta\norm{d^k}
           + O(\norm{d^k}^2)
\end{align*}
where $\zeta>0$ is a value related to the smallest singular value of
$\bm{\omega\g c(x^*)\\ I_{\cA^*}}$. By substituting this lower bound
into the expression for $\rho^k(d^k)$ derived above, and using
\eqref{eq:SLP.d}, we obtain
\begin{equation}\label{eq:rho}
  \rho^k = 1 + \frac{O(\norm{d^k}^2)}{-p^Td^k
              + \omega\norm{\nabla c(x^k)d^k}_1}
        = 1 + O(\norm{d^k})
        = 1 + O(\norm{x^*-x^k}).
\end{equation}

We can now argue recursively to obtain the result.
Suppose that $x^{\bar{k}}$ and $\Delta^{\bar{k}}$ are as defined in
the statement of the theorem, with $\| x^{\bar{k}}-x^* \|$ small
enough that Lemma~\ref{lem:local.conv} holds. By tightening the
requirement on $\| x^{\bar{k}}-x^* \|$ if necessary, we note the following: (i) From
the estimate \eqref{eq:rho}, we have that $\rho^{\bar{k}} \ge \ub{\eta}$, so
that $\Delta^{\bar{k}+1} \ge \Delta^{\bar{k}}$ according to
Algorithm~\ref{alg:SLP}, and $x^{\bar{k}+1} = x^{\bar{k}} +
d^{\bar{k}}$; and (ii) from \eqref{eq:SLP.d}, we have that $\| x^{\bar{k}+1}-x^* \|
= O( \| x^{\bar{k}}-x^* \|^2) \le \|x^{\bar{k}}-x^*\|$, so both
Lemma~\ref{lem:local.conv} and the estimates above continue to hold at
the next iteration $\bar{k}+1$. Thus the recursion continues through
all subsequent iterates $k \ge \bar{k}$. Quadratic convergence follows
from \eqref{eq:SLP.d}, since we have
\[
  \norm{x^{k+1}-x^*} = \norm{x^k+d^k-x^*} = O( \norm{x^k-x^*} ^2).
\]
\qed
\end{proof}


\section{\texorpdfstring{Accelerating \SLP via an Active-Set Heuristic}{Accelerating Sl1LP via an Active-Set Heuristic}} \label{sec:actset}

We now consider the behavior of \SLP as it approaches
``underdetermined'' solutions (see Definition~\ref{def:fully}).  Since
we can expect only linear convergence in these circumstances, we
examine ways to accelerate the method by making use of active-set
estimates and second-order information. In this section, we examine
the elements of this approach in turn. We start with estimation of the
active set in Subsection~\ref{sec:act.set} and of values for the
Lagrange multipliers in
Subsection~\ref{sec:Lagr.mult}. Subsections~\ref{sec:2foundations}
describes the nonlinear system of equations to be solved in the
active-set strategy, while Subsection~\ref{sec:tweak} describes a
heuristic for modifying the active-set estimate if it appears to be
faulty, and concludes with a full specification of the active-set
heuristic.  Subsection~\ref{sec:fitting} describes how the active-set
heuristic is inserted into the \SLP algorithm, and discusses convergence
properties of the enhanced approach.

Our discussion in this section refers both to the original formulation
\eqref{eq:GNO} and the nonsmooth penalty-function form
\eqref{eq:CNSO.constr}. As we have noted, the two forms are equivalent
for sufficiently large choice of penalty parameter $\omega$, and it
makes sense for a locally-convergent phase of the \SLP algorithm to
assume that the iterates have been steered toward a manifold where
$c(x)=0$ holds, so that the linearization of this condition can be
enforced directly (rather than penalized) in computing the steps.

\subsection{Active Set Identification}
\label{sec:act.set}

A vital ingredient of an active-set strategy for \eqref{eq:GNO}
include a reliable means for estimating those bounds that are active
at the solution $x^*$. This issue has been examined in different
contexts; see, for example \cite{Wri93a,HarL04,Lew03a}. It has been
shown that first-order information is often sufficient to make a
reliable identification of the optimal active set, in certain
circumstances.  \cite{ObeW05a} have shown that the optimal active set
of the nonsmooth-penalty formulation of nonlinear programming problem
can be (approximately) identified by linear subproblems under certain
conditions. They show if the current iterate $x^k$ is close enough to
$x^*$ and the trust-region $\Delta^k$ is small enough, then the active
set identified by $x^k+d^k$ is a subset of the optimal active
set. Moreover, if LICQ and strict complementarity conditions hold, and
if the trust-region radius $\Delta^k$ is large enough that $x^*-x^k$
is feasible for the subproblem yet small enough to prevent constraints
inactive at $x^*$ from becoming active in the subproblem, then the
active set identified from the LP subproblem coincides with the
optimal active set.

With this theory in mind, while not attempting to implement it
rigorously, we use the simple heuristic of estimating the active set
from the active set of the subproblem \eqref{eq:SLP.sublp}.
Specifically,
we estimate $\lb{\cA}^*$ and $\ub{\cA}^*$ by the sets $\lb{\cA}_{k+1}$
and $\ub{\cA}_{k+1}$, respectively, where $x^{k+1} = x^k+d^k$ and
\begin{align*}
  \lb{\cA}^{k+1} &:= \set{ i = 1,2,\dotsc,n \given (x^k+d^k)_i =\lb{x}_i }, \\
  \ub{\cA}^{k+1} &:= \set{ i = 1,2,\dotsc,n \given (x^k+d^k)_i =\ub{x}_i },
\end{align*}
as in \eqref{eq:active}.  We observe that on \SLP iterations, there
are fewer and fewer changes to the active set as the iterations
progress. It is therefore reasonable to use these sets as an estimate
of of the optimal active sets once the number of changes drops below a
specified threshold.  Although these estimates may not be exact, they
provide a good starting point for the heuristic discussed in
Subsection~\ref{sec:tweak}, where incremental changes are tried for
the active set in an attempt to restore consistency of the first-order
optimality conditions.


\subsection{Lagrange Multiplier Estimation}
\label{sec:Lagr.mult}

We now discuss how to estimate Lagrange multipliers $\lambda$, $\mu$,
and $\nu$ to initialize the active-set heuristic after an estimate of
the active set becomes available.  Because, as we discuss in
Subsection~\ref{sec:slp.local}, fast local convergence can be obtained
without such estimates in the fully determined case, we consider here
only {\em underdetermined} cases, for which $m+\abs{\cA^*}<n$. We
assume that LICQ holds at $x^*$, that is, the constraint Jacobian
matrix $\bm{\g c(x^*)\\ I_{\cA^*}}$ has full row rank.  The constraint
Jacobian (that is, the basis matrix) for a nondegenerate solution of
the LP subproblem \eqref{eq:SLP.sublp} will contain $n$ rows, with the
extra $n-(m+\abs{\cA^*})$ rows coming from enforcement of additional
bound or trust-region constraints on $d$. We thus make the following
assumption, to ensure that this augmented Jacobian retains its
full-rank property.
\begin{Asn}\label{asn:J.compl}
  All matrices of the form
  \[  \bm{\g c(x^*)\\ I_{\cA^*}\\ E^T} \]
  are nonsingular where $E$ is an $n\times(n-m-\abs{\cA^*})$ matrix
  whose columns are drawn from $\set{e_i \given i \notin \cA^*}$. Thus
  there exists $\xi>0$ such that for all such $E$ we have
  \[ \norm*{\bm{\g c(x^*)\\ I_{\cA^*}\\ E^T}^{-1}} \le \xi. \]
\end{Asn}
Note that this assumption holds trivially for fully determined cases,
since $E$ is null in such cases.

The dual linear program for \eqref{eq:SLP.sublp} is as follows:
\bseq\label{eq:SLP.sublp.d}
\begin{align}
  \omax{\lambda,\mu,\nu}
      & \omega c(x^k)\lambda + (l^k)^T\mu + (u^k)^T\nu \\
  \ost
      & p + \omega\nabla c(x^k)^T\lambda-\mu-\nu
        = 0                                     \label{eq:NLP.sublp.d.constr}\\
      & -\bOne \le \lambda \le \bOne,\quad \mu\ge0,\quad \nu\le 0,
\end{align}
\eseq
where
\[ l^k:=\max ( \lb{x}-x^k,-\Delta^k), \qquad
   u^k:=\min (\ub{x}-x^k,-\Delta^k).
\]
The following theorem states the relation between the solution of
\eqref{eq:SLP.sublp.d} and the optimal Lagrange multipliers
$(\lambda^*,\mu^*,\nu^*)$ under Assumption \ref{asn:J.compl}. This
result rests on an assumption that identification of the active set
via the penalized linear programming model, as described in
\cite{ObeW05a}, has given an accurate result.

\begin{Thm}\label{thm:SLP.lagr}
  Suppose that $(x^*,\lambda^*,\mu^*,\nu^*)$ satisfies
  \eqref{eq:CNSO.KKT} and the optimal active set of \eqref{eq:GNO} is
  identified by $x^k+d^k$. That is, we have
\[
(x^k+d^k)_i = \lb{x}_i \;\;  \Leftrightarrow \;\; i \in \lb{\cA}^*, \qquad
(x^k+d^k)_i = \ub{x}_i \;\;  \Leftrightarrow \;\; i \in \ub{\cA}^*.
\]
If the LP subproblem \eqref{eq:SLP.sublp} has a solution
$(d,\alpha)=(d^k,0)$ which is not degenerate and the Assumption
\ref{asn:J.compl} holds, then provided that $x^k$ is sufficiently
close to $x^*$, the dual problem \eqref{eq:SLP.sublp.d} has a solution
$(\lambda^k,\mu^k,\nu^k)$ such that
\[
(\lambda^k,\mu^k,\nu^k) = (\lambda^*,\mu^*,\nu^*)+O(\norm{x^*-x^k}).
\]
\end{Thm}
\begin{proof}
  Let $\lb{\cA}_\Delta^{k+1}$ and $\ub{\cA}_\Delta^{k+1}$ be defined by
  \begin{align*}
    \lb{\cA}_\Delta^{k+1} &= \set{i=1,2,\dotsc,n \given (x^k+d^k)_i = \lb{x}_i}
                         \cup \set{i=1,2,\dotsc,n \given d_i^k = -\Delta^k},\\
    \ub{\cA}_\Delta^{k+1} &= \set{i=1,2,\dotsc,n \given (x^k+d^k)_i = \ub{x}_i}
                         \cup \set{i=1,2,\dotsc,n \given d_i^k = \Delta^k}.
  \end{align*}
  By the assumptions of the theorem, we have
  $\lb{\cA}^*\subset\lb{\cA}_\Delta^{k+1}$ and
  $\ub{\cA}^*\subset\ub{\cA}_\Delta^{k+1}$. From linear programming
  duality, we have
  \[
    \mu_i^k = 0  \tfor i \notin \lb{\cA}_\Delta^{k+1}, \quad\quad
    \nu_i^k = 0  \tfor i \notin \ub{\cA}_\Delta^{k+1}.
  \]
  We know further that
  \[
    \mu_i^* = 0  \tfor i \in \lb{\cA}_\Delta^{k+1} \backslash \lb{\cA}^*,
\quad\quad
    \nu_i^* = 0  \tfor i \in \ub{\cA}_\Delta^{k+1} \backslash \ub{\cA}^*.
  \]
  Then from \eqref{eq:CNSO.KKT.lagr} and \eqref{eq:NLP.sublp.d.constr},
we have
  \begin{align*}
     \bm{\nabla c(x^k)^T &
         I_{\lb{\cA}_\Delta^{k+1}}^T &
         I_{\ub{\cA}_\Delta^{k+1}}^T}
     \bm{ \omega\lambda^k\\
         -\mu_{\lb{\cA}_\Delta^{k+1}}^k\\
         -\nu_{\ub{\cA}_\Delta^{k+1}}^k}
    &= -p\\
    &\hspace{-70px}= \bm{\nabla c(x^*)^T &
                         I_{\lb{\cA}^*}^T &
                         I_{\ub{\cA}^*}^T}
                     \bm{ \omega\lambda^*\\
                         -\mu_{\lb{\cA}^*}^*\\
                         -\nu_{\ub{\cA}^*}^*} \\
    &\hspace{-70px}= \bm{\nabla c(x^*)^T &
                         I_{\lb{\cA}_\Delta^{k+1}}^T &
                         I_{\ub{\cA}_\Delta^{k+1}}^T}
                     \bm{ \omega\lambda^*\\
                         -\mu_{\lb{\cA}_\Delta^{k+1}}^*\\
                         -\nu_{\ub{\cA}_\Delta^{k+1}}^*}\\
    &\hspace{-70px}= \bm{\nabla c(x^k)^T &
                         I_{\lb{\cA}_\Delta^{k+1}}^T &
                         I_{\ub{\cA}_\Delta^{k+1}}^T}
                     \bm{ \omega\lambda^*\\
                         -\mu_{\lb{\cA}_\Delta^{k+1}}^*\\
                         -\nu_{\ub{\cA}_\Delta^{k+1}}^*}
                     + O(\norm{x^k-x^*}).
  \end{align*}
  Since the solution $(d^k, 0)$ is not degenerate, and by smoothness
  of $c$, we have from Assumption~\ref{asn:J.compl} that the matrix
  $\bm{\nabla
    c(x^k)\\ I_{\lb{\cA}_\Delta^{k+1}\cup\ub{\cA}_\Delta^{k+1}}}^T$ is
  a square, uniformly nonsingular matrix for $x^k$ sufficiently close
  to $x^*$, and thus
  \[
    \bm{ \omega\lambda^k\\
        -\mu_{\lb{\cA}_\Delta^{k+1}}^k\\
        -\nu_{\ub{\cA}_\Delta^{k+1}}^k}
    = \bm{ \omega\lambda^*\\
          -\mu_{\lb{\cA}_\Delta^{k+1}}^*\\
          -\nu_{\ub{\cA}_\Delta^{k+1}}^*}
      + O(\norm{x^k-x^*}).
  \]
  Since $\mu_i^* = \mu_i^k=0$ for $i \notin \lb{\cA}_\Delta^{k+1}$, and $\nu_i^* =
  \nu_i^k = 0$ for $i \notin \ub{\cA}_\Delta^{k+1}$, we have
  \[
    \bm{\lambda^k\\
        \mu^k\\
        \nu^k}
    = \bm{\lambda^*\\
          \mu^*\\
          \nu^*}
      + O(\norm{x^*-x^k}),
  \]
  as desired. \qed
\end{proof}

Theorem \ref{thm:SLP.lagr} indicates that the dual solution of the LP
subproblem can be used as a good approximation of the optimal Lagrange
multipliers of the original CNSO problem \eqref{eq:CNSO.constr}, for
$x^k$ close to $x^*$.  However, the estimated Lagrange multipliers may
violate complementarity conditions in the optimality conditions
\eqref{eq:CNSO.KKT}, due to the presence of trust-region
constraints. To fix these violations, we drop the values of dual
solution corresponding to the trust-region constraints and define
$(\lambda^{k+1}, \mu^{k+1}, \nu^{k+1})$ to start the active-set  heuristic
at iterate $x^{k+1} = x^k+d^k$ as below:
\begin{equation}\label{eq:lag.mult.approx}
  \lambda^{k+1} = \lambda^k, \quad
  \mu^{k+1}_i   = \begin{cases}
                     \mu^k_i & i \in \lb{\cA}^{k+1}, \\
                     0       & \tow,
                   \end{cases} \quad
  \nu^{k+1}_i   = \begin{cases}
                     \nu^k_i & i \in \ub{\cA}^{k+1}, \\
                     0       & \tow.
                   \end{cases}
\end{equation}


\subsection{The Active-Set Heuristic} \label{sec:2foundations}

The active-set heuristic consists of a sequence of Newton step on a
system of nonlinear equations like \eqref{eq:KKT.eq}, with
$\lb{\cA}^*$ and $\ub{\cA}^*$ replaced by current estimates of these
active sets, and the Newton linearization taking place around the
latest primal-dual iterate
$(\tilde{x},\tilde{\lambda},\tilde{\mu},\tilde{\nu})$, which is
initialized to the primal-dual point
$(x^k+d^k,\lambda^{k+1},\mu^{k+1},\nu^{k+1})$ obtained as in the
previous subsection. We initialize the active set estimates
$\lb{\cA}'$ and $\ub{\cA}'$ to $\lb{\cA}^{k+1}$ and $\ub{\cA}^{k+1}$,
respectively, and define
\[
\lb{\cI}' = \set{1,2,\dotsc,n} \backslash \lb{\cA}', \quad
\ub{\cI}' = \set{1,2,\dotsc,n} \backslash \ub{\cA}', \quad
\cI' = \set{1,2,\dotsc,n} \backslash (\lb{\cA}' \cup \ub{\cA}').
\]
as in \eqref{eq:inactive}.  The Newton equations are as follows:
\begin{equation}\label{eq:KKT.LS}
  \bm{\tilde{H}                  & \tilde{A} & -I          & -I          \\
      \tilde{A}^T              &     &             &             \\
      I_{\lb{\cA}',\cdot} &     &             &             \\
      I_{\ub{\cA}',\cdot} &     &             &             \\
                           &     & I_{\lb{\cI}',\cdot} &             \\
                           &     &             & I_{\ub{\cI}',\cdot}}
  \bm{\Delta x\\
      \Delta \lambda\\
      \Delta \mu\\
      \Delta \nu}
  = -\bm{\nabla_x\cL(\tilde{x},\tilde{\lambda},\tilde{\mu},\tilde{\nu})\\
         c(\tilde{x})\\
         (\tilde{x}-\lb{x})_{\lb{\cA}'}\\
         (\tilde{x}-\ub{x})_{\ub{\cA}'}\\
         \tilde{\mu}_{\lb{\cI}'}\\
         \tilde{\nu}_{\ub{\cI}'}}
\end{equation}
where
\[
  \tilde{H}    = \sum_{i=1}^m\tilde{\lambda}_i\nabla^2 c_i(\tilde{x}), \quad
 \tilde{A}           = \nabla c(\tilde{x})^T.
\]
(Note that the coefficient matrix in \eqref{eq:KKT.LS} is square.)
Denoting the solution of \eqref{eq:KKT.LS} by $(\Delta x', \Delta
\lambda', \Delta \mu', \Delta \nu')$, we define a provisional estimate
of the next iterate by
\begin{equation}\label{eq:KKT.soln.update}
  (x',\lambda',\mu',\nu')
  = (\tilde{x}, \tilde{\lambda}, \tilde{\nu}, \tilde{\mu})
    + (\Delta x', \Delta \lambda', \Delta\mu', \Delta\nu').
\end{equation}
We accept $(x',\lambda',\mu',\nu')$ as the new iterate only if the
following set of inequalities (analogous to \eqref{eq:KKT.ineq}) hold:
\bseq \label{eq:KKT.ineq.k}
\begin{align}
  \lb{x}_i < x_i' < \ub{x}_i & \tfor i\in\cI', \label{eq:KKT.ineq.k.x} \\
  \mu_i' \ge 0               & \tfor i\in\lb{\cA}', \label{eq:KKT.ineq.k.mu} \\
  \nu_i' \le 0               & \tfor i\in\ub{\cA}'. \label{eq:KKT.ineq.k.nu}
\end{align}
\eseq
We require in addition that the new iterate improves the norm of the
algebraic KKT conditions \eqref{eq:kkt.L}, \eqref{eq:kkt.c} by at
least a factor of $1/2$, that is,
\begin{equation} \label{eq:kkt.impr}
\norm*{ \bm{ \g_x \cL(x',\lambda',\mu',\nu') \\ c(x') } } \le \frac12
\norm*{ \bm{ \g_x \cL (\tilde{x}, \tilde{\lambda}, \tilde{\nu}, \tilde{\mu}) \\
c(\tilde{x}) } }.
\end{equation}
If both these conditions hold, we replace
$(\tilde{x},\tilde{\lambda},\tilde{\mu},\tilde{\nu})$ by
$(x',\lambda',\mu',\nu')$, leave the active-set estimates $\lb{\cA}'$
and $\ub{\cA}'$ unchanged, and take another active-set step.

If condition \eqref{eq:KKT.ineq.k} fails to hold, we try to ``tweak''
the active set and re-solve the Newton equations for the new active
sets, as described in the next subsection.

Note that the system \eqref{eq:KKT.LS} contains a great deal of
structure. By performing various block eliminations and substitutions,
we can reduce it significantly in size. Some of the variables can be
obtained directly, as follows
\begin{equation}\label{eq:KKT.soln.pre}
\begin{aligned}
  \Delta x_{\lb{\cA}'}'   &= (\lb{x}-\tilde{x})_{\lb{\cA}'}, \\
  \Delta x_{\ub{\cA}'}'   &= (\ub{x}-\tilde{x})_{\ub{\cA}'}, \\
  \Delta \mu_{\lb{\cI}'}' &= -\mu_{\lb{\cI}'}, \\
  \Delta \nu_{\ub{\cI}'}' &= -\nu_{\ub{\cI}'}.
\end{aligned}
\end{equation}
There are more variables that appear in just one equation, whose
values can be defined as follows:
\begin{equation}\label{eq:KKT.soln.post}
\begin{aligned}
  \Delta \mu_{\lb{\cA}'}'
    &= (\tilde{H}\Delta x' + \tilde{A} \Delta \lambda' + \nabla_x
\cL(\tilde{x}, \tilde{\lambda}, \tilde{\nu}, \tilde{\mu}))_{\lb{\cA}'},\\
  \Delta \nu_{\ub{\cA}'}'
    &= (\tilde{H}\Delta x' + \tilde{A} \Delta \lambda' + \nabla_x
\cL(\tilde{x}, \tilde{\lambda}, \tilde{\nu}, \tilde{\mu}))_{\ub{\cA}'}.
\end{aligned}
\end{equation}
The remaining variables can be obtained by solving the following
linear system:
\begin{equation}\label{eq:KKT.LS.reduced}
  \bm{\tilde{H}_{\cI',\cI'}     & \tilde{A}_{\cI',\cdot} \\
      (\tilde{A}_{\cI',\cdot})^T &    0          }
  \bm{\Delta x_{\cI'}'\\
      \Delta\lambda'}
  = -
  \bm{[\nabla_x \cL(\tilde{x}, \tilde{\lambda}, \tilde{\nu}, \tilde{\mu})]_{\cI'}\\
      c(\tilde{x})}.
\end{equation}
The coefficient matrix in \eqref{eq:KKT.LS.reduced} is symmetric
indefinite, so we can perform an factorization involving a lower
triangular matrix $L$ and a block diagonal matrix $D$:
\begin{equation} \label{eq:ldl}
  LDL^T = \bm{\tilde{H}_{\cI',\cI'}     & \tilde{A}_{\cI',\cdot} \\
      (\tilde{A}_{\cI',\cdot})^T &     0         }.
\end{equation}
Having calculated this factorization, the solution of
\eqref{eq:KKT.LS.reduced} can be computed by performing two triangular
substitutions and some other simple operations.


\subsection{Active Set Adjustment}
\label{sec:tweak}

In Subsection~\ref{sec:act.set}, we discussed taking the active set
from the latest LP subproblem as our estimate of the optimal active
set. We may find that for the step produced in \eqref{eq:KKT.LS},
\eqref{eq:KKT.soln.update} does not satisfy the inequality conditions
\eqref{eq:KKT.ineq.k}, which may be an indication that the current
active set estimate is not optimal. Since it is probably not too
different from the optimum, however, we propose another heuristic for
making a few changes to it, rather than discarding the step completely
and returning to the \SLP algorithm.

The following rules are used to modify the active-set estimates
$\lb{\cA}'$ and $\ub{\cA}'$.
\begin{itemize}
\item If $x'_i<\lb{x}_i$, we add $i$ to $\lb{\cA}'$; if
  $x'_i>\ub{x}_i$, we add $i$ to $\ub{\cA}'$.
\item If $\mu'_i<0$, we remove $i$ from $\lb{\cA}'$.
\item If $\nu'_i>0$, we remove $i$ from $\ub{\cA}'$.
\end{itemize}
Following these changes, we calculate the Newton step from
\eqref{eq:KKT.LS}, with the modified sets $\lb{\cA}'$ and $\ub{\cA}'$,
recalculate the provisional new iterate $(x',\lambda',\mu',\nu')$ from
\eqref{eq:KKT.soln.update}, and retest the inequality conditions
\eqref{eq:KKT.ineq.k}.  We declare ``success'' when the inequality
conditions are satisfied and the sufficient decrease test
\eqref{eq:kkt.impr} holds. Otherwise, we repeat the adjustment
procedure up to a predefined number of iterations. We declare
``failure'' and return to \SLP iterations if
\begin{itemize}
\item the conditions \eqref{eq:KKT.ineq.k} are not satisfied at any
  of these steps,
\item the conditions \eqref{eq:KKT.ineq.k} are satisfied, but the
  improvement by the estimate $(\tilde{x}, \tilde{\lambda},
  \tilde{\mu}, \tilde{\nu})$ is insufficient. (i.e. \eqref{eq:kkt.impr}
  is not satisfied.)
\item we re-encounter estimates $\lb{\cA}'$ and $\ub{\cA}'$ that
  had been tried already at one of the previous steps, or
\item the total number of components in the active set ($|
  \lb{\cA}'| + |\ub{\cA}'|$) grows larger than $n-m$.
\end{itemize}
The second condition is invoked to prevent cycling among a set of
choices for the active-set estimates.  When the third condition holds,
the solution is ``overdetermined'' by its constraints, and the
coefficient matrix in \eqref{eq:KKT.LS} becomes structurally singular.

Rather that solving \eqref{eq:KKT.LS} from scratch after modification
of the active sets, we re-use the factorization \eqref{eq:ldl},
modifying the system \eqref{eq:KKT.LS.reduced} by adding columns to
account for the changes to the active-set estimates.
Any revised problem with the altered active set approximation can be
written in the following form:
\begin{equation}\label{eq:KKT.LS.AUG}
  \left[\begin{array}{cc|c}
    \tilde{H}_{\cI',\cI'}     & \tilde{A}_{\cI',\cdot} & \mr{2}{*}{V} \\
    (\tilde{A}_{\cI',\cdot})^T &      0       &  \\
    \hline
    \multicolumn{2}{c|}{X^T}  & S \\
  \end{array}\right]
  \left[\begin{array}{c}
    \Delta x_{\cI'}\\ \Delta\lambda \\\hline \cdot
  \end{array}\right]
  = - \left[\begin{array}{c}
        [\nabla_x \cL(\tilde{x},\tilde{\lambda},\tilde{\mu},\tilde{\nu})]_{\cI'} \\ \c(\tilde{x}) \\ \hline s
      \end{array}\right]
\end{equation}
where the matrices $X$, $V$, $S$ and the vector $s$ capture the
changes to the active set.
Below, we describe the five possible ways in which the active set can
be altered, and show how to define the system \eqref{eq:KKT.LS.AUG} to
account for these cases.

\begin{enumerate}[label=(\emph{\alph*})]
\item \label{enum:tweak.1} When $x_i$, for some $i \notin \cA'$, moves
  from the interior of its box constraint to beyond the lower bound.
  Since the constraint $x_i\ge \lb{x}_i$ needs to be active, the value
  of $\Delta x_i$ should be $\lb{x}_i-x_i^k$. Also the corresponding
  Lagrange multiplier $\mu_i$ is allowed to move away from zero.  The
  modified linear system can thus be defined as follows:
  \[
    \left[\begin{array}{cc|c}
      \tilde{H}_{\cI',\cI'}     & \tilde{A}_{\cI',\cdot} & (e_i)_{\cI'} \\
      (\tilde{A}_{\cI',\cdot})^T &      0       &  0 \\
      \hline
      (e_i)_{\cI'}^T &  0 &  0\\
    \end{array}\right]
    \left[\begin{array}{c}
      \Delta x_{\cI'}\\
      \Delta\lambda  \\\hline
      -\Delta\mu_i
    \end{array}\right]
    = - \left[\begin{array}{c}
          [\nabla_x \cL(\tilde{x},\tilde{\lambda},\tilde{\mu},\tilde{\nu})]_{\cI'} \\
          c(\tilde{x}) \\\hline
          \tilde{x}_i-\lb{x}_i
        \end{array}\right].
  \]
\item \label{enum:tweak.2} When $x_i$, for some $i \notin \cA'$, moves
  from the interior of its box constraint to the upper bound, a similar
  construction yields the following modified system:
  \[
    \left[\begin{array}{cc|c}
      \tilde{H}_{\cI',\cI'}     & \tilde{A}_{\cI',\cdot} & (e_i)_{\cI'} \\
      (\tilde{A}_{\cI',\cdot})^T &    0         & 0  \\
      \hline
      (e_i)_{\cI'}^T &  0 &  0 \\
    \end{array}\right]
    \left[\begin{array}{c}
      \Delta x_{\cI'}\\
      \Delta\lambda \\\hline
      -\Delta \nu_i
    \end{array}\right]
    = - \left[\begin{array}{c}
          [\nabla_x \cL(\tilde{x},\tilde{\lambda},\tilde{\mu},\tilde{\nu})]_{\cI'} \\
          c(\tilde{x}) \\\hline
          \tilde{x}_i-\ub{x}_i
        \end{array}\right].
  \]
\item \label{enum:tweak.3} When, for some $i \in \cA'$, a component
  $\nu_i$ or $\mu_i$ moves from one nonzero value to another nonzero
  value with a different sign, we take it as an indication that the
  corresponding component of $x$ should move away from its bound. We
  thus set the Lagrange multiplier component in question to zero, and
  allow $\Delta x_i$ to become nonzero, The augmented linear system
  has the following form:
  \[
    \left[\begin{array}{cc|c}
      \tilde{H}_{\cI',\cI'}     & \tilde{A}_{\cI',\cdot} & \tilde{H}_{\cI',i} \\
      (\tilde{A}_{\cI',\cdot})^T &      0       & (\tilde{A}_{i,\cdot})^T \\
      \hline
      \tilde{H}_{i,\cI'} & \tilde{A}_{i,\cdot}  & \tilde{H}_{i,i} \\
    \end{array}\right]
    \left[\begin{array}{c}
      \Delta x_{\cI'}\\
      \Delta\lambda \\\hline
      \Delta x_i
    \end{array}\right]
    = - \left[\begin{array}{c}
          [\nabla_x \cL(\tilde{x},\tilde{\lambda},\tilde{\mu},\tilde{\nu})]_{\cI'} \\
          c(\tilde{x}) \\ \hline
          [\nabla_x \cL(\tilde{x},\tilde{\lambda},\tilde{\mu},\tilde{\nu})]_i
        \end{array}\right].
  \]
If we permute the rows and columns of this system, we can recover a
structure similar to \eqref{eq:KKT.LS}.

\item \label{enum:tweak.4} When, for some $i \in \ub{\cA}'$, $x_i$
  moves from its upper bound to its lower bound, we redefine $\Delta
  x_i = \lb{x}_i - \ub{x}_i$, and define the augmented system as
  follows:
  \[
    \left[\begin{array}{cc|c}
      \tilde{H}_{\cI',\cI'}     & \tilde{A}_{\cI',\cdot}  & \tilde{H}_{\cI',i} \\
      (\tilde{A}_{\cI',\cdot})^T &       0      & (\tilde{A}_{i,\cdot})^T \\
      \hline
      0  & 0  & 1 \\
    \end{array}\right]
    \left[\begin{array}{c}
      \Delta x_{\cI'}\\
      \Delta\lambda \\\hline
      \Delta x_i
    \end{array}\right]
    = - \left[\begin{array}{c}
           [\nabla_x \cL(\tilde{x},\tilde{\lambda},\tilde{\mu},\tilde{\nu})]_{\cI'}  \\
          c(\tilde{x}) \\\hline
          \ub{x}_i-\lb{x}_i
        \end{array}\right].
  \]
\item \label{enum:tweak.5} When, for some $i \in \lb{\cA}'$, $x_i$
  moves from its lower bound to its upper bound, we modify the
augmented system as follows:
  \[
    \left[\begin{array}{cc|c}
      \tilde{H}_{\cI',\cI'}     & \tilde{A}_{\cI',\cdot} & \tilde{H}_{\cI',i} \\
      (\tilde{A}_{\cI',\cdot})^T &    0         & (\tilde{A}_{i,\cdot})^T \\
      \hline
      0  & 0  & 1 \\
    \end{array}\right]
    \left[\begin{array}{c}
      \Delta x_{\cI'}\\
      \Delta\lambda \\\hline
      \Delta x_i
    \end{array}\right]
    = - \left[\begin{array}{c}
          [\nabla_x \cL(\tilde{x},\tilde{\lambda},\tilde{\mu},\tilde{\nu})]_{\cI'} \\
          c(\tilde{x}) \\\hline
          \lb{x}_i-\ub{x}_i
        \end{array}\right].
  \]
\end{enumerate}

These modifications can be combined when multiple changes to the active
set are made on a single step. Each such change results in one extra
row and column being added to the reduced augmented system.
Appendix~\ref{appdx:AugLS} contains further details on how the
factorization \eqref{eq:ldl} for the original coefficient matrix can
be leveraged to solve the systems above efficiently.

The full specification of the active-set heuristic appears as
Algorithm~\ref{alg:asi}.

\begin{algorithm}
\caption{Active-Set Heuristic} \label{alg:asi}
\begin{algorithmic}[1]

\Require

\Statex Primal-dual iterate
$(\tilde{x},\tilde{\lambda},\tilde{\mu},\tilde{\nu})$ and active set
estimates $\lb{\cA}'$ and $\ub{\cA}'$;
\Statex Maximum number of iterations for tweaking process: $T$;

\Ensure

\Statex {\bf Either} declare failure and return to iterating in
  Algorithm~\ref{alg:SLP}, \par
  {\bf or} produce a primal-dual solution
  $(x^*,\lambda^*,\mu^*,\nu^*)$;

\medskip

\Loop
\State Calculate candidate iterate  $(x',\lambda',\mu',\nu')$ from
  \eqref{eq:KKT.LS} and \eqref{eq:KKT.soln.update};

\State {TweakCounter $\gets$ 0;}
\While {$(x',\lambda',\mu',\nu')$ does {\em not} satisfy conditions \eqref{eq:KKT.ineq.k}}
  \State {TweakCounter $\gets$ TweakCounter + 1;}
  \If {TweakCounter $> T$}
    \State {\bf stop} and return to Algorithm~\ref{alg:SLP};
      \Comment{No suitable active set found.}
  \EndIf
  \State {Adjust active set estimates $\lb{\cA}'$ and $\ub{\cA}'$ using the rules described in Subsection~\ref{sec:tweak};}
  \If{$\lb{\cA}'$ and $\ub{\cA}'$ were encountered on
    a previous iteration of this while loop \par
    \hskip\algorithmicindent {\bf or} $| \lb{\cA}'| + |\ub{\cA}'| > n-m$
    }
  \State {\bf stop} and return to Algorithm~\ref{alg:SLP};
    \Comment{No suitable active set found.}
  \EndIf
  \State Calculate candidate iterate $(x',\lambda',\mu',\nu')$ for the
  adjusted sets $\lb{\cA}'$ and $\ub{\cA}'$;
\EndWhile

\If{\eqref{eq:kkt.impr} does {\em not} hold}
  \State {\bf stop} and return to Algorithm~\ref{alg:SLP};
    \Comment{Insufficient improvement.}
\EndIf

\State Set $(\tilde{x},\tilde{\lambda},\tilde{\mu},\tilde{\nu}) \leftarrow
(x',\lambda',\mu',\nu')$;
\If{
$\norm*{ \bm{ \g_x \cL(\tilde{x},\tilde{\lambda},\tilde{\mu},\tilde{\nu})  \\ c(\tilde{x}) } } \le \mbox{\bf tol}
$
}
  \State {$(x^*,\lambda^*,\mu^*,\nu^*) \gets (\tilde{x},\tilde{\lambda},\tilde{\mu},\tilde{\nu})$;} \Comment{Success!}
  \State{\bf stop};
\EndIf
\EndLoop

\end{algorithmic}
\end{algorithm}

\subsection{Inserting the Active-Set Heuristic into the \SLP Algorithm}
\label{sec:fitting}

The active-set heuristic can be invoked after an iteration of the \SLP
algorithm is completed, if we decide that the active sets $\lb{\cA}^k$
and $\ub{\cA}^k$ have settled down (for example, they have changed by
fewer than ten components). To invoke the active-set heuristic after a
successful iteration $k$ of \SLP, we estimate Lagrange multipliers
according to the procedure of Subsection~\ref{sec:Lagr.mult}, then
invoke Algorithm~\ref{alg:asi} with
$(\tilde{x},\tilde{\lambda},\tilde{\mu},\tilde{\nu}) =
(x^k+d^k,\lambda^{k+1},\mu^{k+1},\nu^{k+1})$, $\lb{\cA}' =
\lb{\cA}^{k+1}$, and $\ub{\cA}' = \ub{\cA}^{k+1}$.

If Algorithm~\ref{alg:asi} exits with a failure of the active-set
heuristic, we simply return to the \SLP algorithm and pick up where we
left off, with the latest values of $x^k$ and $\Delta^k$, and proceed
with further iterations according to Algorithm~\ref{alg:SLP}. Thus, it
is clear that the global convergence theory developed in
Section~\ref{sec:slp} continues to hold.

Otherwise, if we remain in the active-set heuristic indefinitely, the
conditions \eqref{eq:KKT.ineq.k} together with the decrease condition
\eqref{eq:kkt.impr} ensure that any limit point of this heuristic
satisfies the first-order optimality conditions \eqref{eq:kkt} of
\eqref{eq:GNO}. Although we refrain from offering a formal theory for
the rate of convergence, it is clear that if the active-set heuristic is
invoked in the neighborhood of a point $x^*$ at which first-order
conditions, LICQ, strict complementarity, and second-order sufficient
conditions for \eqref{eq:GNO} are satisfied, it can be expected to
converge quadratically, provided that $c$ is twice continuously
differentiable near $x^*$.

\section{Feasibility Restoration for Power Systems}
\label{sec:formulation}

In this section we outline the problem of restoring feasibility to a
power network by optimal load shedding. The AC power flow model is
defined in Subsection~\ref{sec:acmodel}, while we formulate the
optimal load-shedding problem as a problem of the form \eqref{eq:GNO}
in Subsection~\ref{sec:feasrest}.

We make use throughout this section of the following notations:
\begin{itemize}
  \item Set of buses: $\cN$
  \item Set of generators: $\cG \subseteq \cN$
  \item Set of demand buses: $\cD \subseteq \cN$
  \item Index of the slack (reference) bus: $s\in\cN$
  \item Set of lines: $\cL \subseteq \cN\times\cN$
  \item Unit imaginary number: $j$
  \item Complex power at bus $i\in\cN$: $P_i+jQ_i$. ($P_i$ is the
    active power and $Q_i$ the reactive power).
  \item Complex voltage at bus $i\in\cN$: $V_ie^{j\theta_i}$. ($V_i$
    is the voltage magnitude and $\theta_i$ is the voltage angle.)
  \item Difference of $\theta_i$ and $\theta_k$: $\theta_{ik}=\theta_i-\theta_k$
  \item Admittance for line $(i,k)$ (that is, the $(i,k)$ element of the
line admittance matrix): $G_{ik}+jB_{ik}$.
\end{itemize}
We note that $\set{ \cG, \cD, \set{s} }$ is a partition of $\cN$, i.e.
the sets $\cG, \cD$ and $\set{s}$ are mutually disjoint and $\cD \cup
\cG \cup \set{s} = \cN$.

\subsection{Power Flow Problems}
\label{sec:acmodel}

To operate a power system we need to know the voltages and powers at
all buses, captured in the vector $(V,\theta,P,Q)$, that satisfy the
following power-flow balance equations arising from Kirchhoff's laws:
\begin{equation} \label{eq:acpf}
  \bm{ F^P(V,\theta)\\
       F^Q(V,\theta)  } = 0,
\end{equation}
where the $i$th entries of $F^P$ and $F^Q$ are nonlinear functions defined as
\bseq \label{eq:ac}
\begin{align}
\label{eq:ac.p}
 F^P_i(V,\theta) &:= V_i\sum_{\mathclap{k:(i,k)\in\cL}}V_k
            (G_{ik}\cos{(\theta_{ik})}+B_{ik}\sin{(\theta_{ik})})-P_i\\
\label{eq:ac.q}
 F^Q_i(V,\theta) &:= V_i\sum_{\mathclap{k:(i,k)\in\cL}}V_k
            (G_{ik}\sin{(\theta_{ik})}-B_{ik}\cos{(\theta_{ik})})-Q_i.
\end{align}
\eseq
The power flow problem is to find a solution $(V,\theta,P,Q)$ to the
above equations, where two out of each quadruplet
$(V_i,\theta_i,P_i,Q_i)$ are known for each bus. The known quantities
are varied by the type of bus. For generators $i \in \cG$, voltage
magnitude $V_i$ and the active power $P_i$ that the generators can
produce are known. For the slack bus $s$, the voltage magnitude $V_s$
and angle $\theta_s$ are known. For load buses $i \in \cD$, the active
power $P_i$, and reactive power $Q_i$ are known. Since $Q_i$ for $i
\in \cG$ can be determined directly from \eqref{eq:ac.q}, we can
reduce \eqref{eq:acpf} to a system of $2\abs{\cD}+\abs{\cG}$ nonlinear
equations with $2\abs{\cD}+\abs{\cG}$ variables ($V_{\cD},
\theta_{\cG\cup\cD}$):
\begin{equation}\label{eq:intro.pflow}
F(V,\theta) =
\bm{F_\cG^P(V,\theta)\\ F_\cD^P(V,\theta)\\ F_\cD^Q(V,\theta)} = 0,
\end{equation}
where $V_s$, $\theta_s$, $V_{\cG}$, $P_{\cG}$, $P_{\cD}$ and $Q_{\cD}$
are given.

Newton's method is widely used for solving \eqref{eq:intro.pflow}. If
a problem is well-conditioned and a sufficiently good initial starting
point is given, it converges to a solution in a few
iterations. Moreover, since first derivatives can be computed easily
for this system, and the Jacobian is quite sparse, it can be
implemented efficiently.

More robust methods to solve ill-conditioned or badly-initialized
power flow problems have been studied by several authors
\citep{IwaT81,TriPM82,Mil09}, including damped Newton methods and
homotopy methods.
There is no known Newton-based method that can detect reliably the
unsolvability of a power flow problem, but in practice, failure of
Newton's methods from a reasonable starting point is strong evidence
of nonexistence of a solution. A method based on Semidefinite programming
(SDP), which is studied by \cite{LavL12} for OPF problems, has been proposed
recently as a more rigorous means of finding solutions and proving
nonexistence, but it cannot resolve all cases (see \cite{MolLD12}).

\subsection{Formulation of Feasibility Restoration for Power Systems}
\label{sec:feasrest}

%

For systems in which \eqref{eq:intro.pflow} does not have a solution,
we have to consider adjustments of demand or generation to restore
feasibility of these equations. Previous works \citep{BarS01a, Ove95}
seek a minimal adjustment in the sense of an $\ell_2$-norm (Euclidean
distance), but the resulting solution is unappealing from an
operational perspective, since it requires loads to be adjusted at
many demand nodes. We consider instead adjustments that minimize an
$\ell_1$-norm, which tend to yield adjustments at fewer nodes and
hence to be potentially more appealing in practice. The $\ell_1$-norm
formulation also leads itself well to the linear-programming-based
techniques described in this paper.


While the minimum adjustment of demand and generation approach
captures the basic concept of feasibility restoration, we need to
refine it by taking certain practicalities into account. Power flow
problem formulations often omit restrictions on voltage magnitudes on
buses, since these do not usually enter into consideration under
normal conditions. Since we are dealing with stressed and disrupted
networks here, it makes good sense to include them in our
formulations, to exclude solutions that would not be operational by
restricting the voltages to a certain range.  Another factor to
consider, since we intend the formulation to provide a practical
indication of how feasibility can be restored, is to impose practical
constraints on the adjustments indicated by the model.  For example,
when load-shedding is needed on a demand bus, the amount of the
load-shedding should not exceed the total amount of demand at the bus.
Also, since the active and reactive powers on a bus are closely
related, they should be adjusted by the same fraction.  With these
additional restrictions, we obtain the following formulation for
feasibility restoration problem:
%
\bseq \label{eq:feas.rest}
\begin{alignat}{2}
  \omin{\substack{V_\cD,\theta_{\cD\cup\cG},\\ \sigma_\cG^+,  \sigma_\cG^-, \rho_\cD}}
       & \mathrlap{\sum_{i\in \cG} {\abs{P_i}(\sigma_i^++\sigma_i^-)}
         + \sum_{i\in \cD} {(\abs{P_i}+\abs{Q_i})\rho_i}}\label{eq:feas.rest.obj}\\
  \ost & F^P_i(V,\theta)-\abs{P_i}(\sigma_i^+-\sigma_i^-) = 0
                                 & \quad & i\in \cG\label{eq:feas.rest.pv.p}\\
       & F^P_i(V,\theta)-\abs{P_i}\rho_i = 0
                                      && i\in \cD\label{eq:feas.rest.pq.p}\\
       & F^Q_i(V,\theta)-\abs{Q_i}\rho_i = 0
                                      && i\in \cD\label{eq:feas.rest.pq.q}\\
       & \lb{V} \le V_i \le \ub{V}
                                      && i\in \cD\label{eq:feas.rest.v}\\
       & 0 \le \sigma_i^+ \le \ub{\sigma}_i^+
                                      && i\in \cG\label{eq:feas.rest.sigma.p}\\
       & 0 \le \sigma_i^- \le \ub{\sigma}_i^-
                                      && i\in \cG\label{eq:feas.rest.sigma.m}\\
       & 0 \le \rho_i \le \ub{\rho}_i
                                      && i\in \cD, \label{eq:feas.rest.rho}
\end{alignat}
\eseq
where $\lb{V}$ and $\ub{V}$ are the lower and upper limits on voltage
magnitudes of demand buses, $\ub{\sigma}_i^\pm$, $i \in \cG$ are the
bounds on the active power adjustments of generators and
$\ub{\rho}_i$, $i \in \cD$ are bounds on the active and reactive power
adjustments of demand nodes. We have assumed here that $P_i\le 0$ and
$Q_i\le0$ for $i\in\cD$, and $P_i>0$ for $i\in\cG$, following
convention.
We note the following points.
\bi
\item The problem has the form \eqref{eq:GNO} for which we describe
  and analyze algorithms in earlier sections. For the purpose of
  applying Algorithm~\ref{alg:SLP}, \eqref{eq:feas.rest} must be
  rewritten as \eqref{eq:CNSO.constr}, whose corresponding LP
  subproblem is \eqref{eq:SLP.sublp}.
\item Constraints \eqref{eq:feas.rest.pv.p}, \eqref{eq:feas.rest.pq.p},
  and \eqref{eq:feas.rest.pq.q} represent relaxations of the power
  flow equations \eqref{eq:intro.pflow} in which the loads $P_\cG$,
  $P_\cD$, and $Q_\cD$ are modified by a certain relative amount,
  captured by the variables $\sigma_i^+$, $\sigma_i^-$, and $\rho_i$.
\item Constraint \eqref{eq:feas.rest.pv.p} ensures that power
  generation can be either increased or decreased, but
  \eqref{eq:feas.rest.pq.p} and \eqref{eq:feas.rest.pq.q} ensure that
  loads at demand nodes can only decrease.
\item The same variable $\rho_i$ is used in the active and reactive
  power balance equations \eqref{eq:feas.rest.pq.p}
  and \eqref{eq:feas.rest.pq.q}, since it makes operational sense in
  many situations for active and reactive load shedding to occur in
  the same fraction.
\item Box constraints on the load shedding variables
  \eqref{eq:feas.rest.sigma.p}, \eqref{eq:feas.rest.sigma.m}, and
  \eqref{eq:feas.rest.rho} ensure that adjustments cannot exceed
  user-defined limits.  (Upper bounds $\ub{\sigma}_i^+$,
  $\ub{\sigma}_i^-$, and $\ub{\rho}_i$ should not exceed 1.)
\item The bounds \eqref{eq:feas.rest.v} guarantee that voltage levels
  are operationally viable in the solution.
\item In the objective \eqref{eq:feas.rest.obj}, we weight the
  coefficient of $\rho_i$ with the sum of original active power
  demand and reactive power demand $(|P_i|+|Q_i|)$.
\ei

The feasibility restoration problem \eqref{eq:feas.rest} is neither
linear nor convex, so we can guarantee only a local solution. The
problem generalizes \eqref{eq:intro.pflow} in that if a solution of
the latter problem exists, it will yield a global solution of
\eqref{eq:feas.rest} with an objective of zero when we set
$\sigma^+_i=\sigma^-_i=0$ for $i \in \cG$ and $\rho_i=0$ for $i \in
\cD$, provided the voltage constraints \eqref{eq:feas.rest.v} are
satisfied.  Moreover, by the well-known sparsity property induced by
the $\ell_1$ objective, we expect few of the components of
$\sigma_{\cG}^+$, $\sigma_{\cG}^-$, and $\rho_\cD$ to be nonzero at a
typical solution of \eqref{eq:feas.rest}.

We note that the formulation \eqref{eq:feas.rest} can be enhanced (at
the cost of some additional complexity in the model) by adding limits
on current flows on the lines, which can be modeled by a combination
of equality constraints and nonnegative slack variables.

The solution of \eqref{eq:feas.rest} suggests to the grid operator a
feasible operating point and a load-shedding pattern that can be used
to attain this point.

\section{Experimental Results} \label{sec:result}

In this section we present the experimental results for our \SLP
algorithm applied to the CNSO reformulation of feasibility restoration
problem \eqref{eq:feas.rest}. The algorithm is implemented in
\MATLAB\footnote{Version 8.1.0.604 (R2013a)} on MacBook Pro (2.6 GHz
Intel Core i7 with 8GB RAM) with \CPLEX\footnote{Version 12.6} as a
linear programming solver. Specifically, the \CPLEX Class API for
\MATLAB is used to exploit the warm start feature of LP solver. The
\MATLAB functions provided by \MATPOWER\footnote{Version 4.1}
\citep{ZimMT11} are used to read power systems data and to compute the
values required to formulate the problems (such as the admittance
matrix and the derivatives of the power flow equations).  The data
sets included in \MATPOWER package, which contains the power systems
data from Power Systems Test Case Archive~\citep{PSTCA}, are used for
the experiments, but they are modified to create infeasible instances.
(The original data is feasible for the power flow problem, so leads to
a zero feasible objective in \eqref{eq:feas.rest} and its CNSO
reformulation.) Our implementation includes pre-compiled code (\CPLEX)
and \MATLAB code.  We compare it with \IPOPT\footnote{Version 3.11.7}
\citep{WacB06}, an interior-point solver written in \verb!C++!  and
compiled for \MATLAB with the MA27 linear solver; and MIPS, a
\MATPOWER Interior Point Solver written in \MATLAB code.

For the \SLP phase, the following values are used for the parameters
that appear in Algorithm~\ref{alg:SLP} and \eqref{eq:SLP.converged}:
\[
\begin{matrix}
  \lb{\eta} = 0.5, \quad
  \ub{\eta} = 0.25, \quad
  \lb{\rho} = 0.1, \quad
  \epsilon  = 10^{-3},\\
  c_1 = 2, \quad
  c_2 = 0.5, \quad
  \Delta^0 = \ub{\Delta} = 1, \quad
  \lb{\Delta} = 10^{-5}.
\end{matrix}
\]
The penalty parameter $\omega$ is set to
$\max\{10\abs{P_{\cD\cup\cG}}, 10\abs{Q_\cG}\}$, which is chosen to
ensure that all equality constraints are satisfied, the tolerance to
check the activities and violations of constraints is set to $10^{-9}$
for all experiments. When the active-set heuristic is used, we set
$T=10$ for the maximum number of tweaking iterations. The choice is
based on the fact that we want to terminate the tweaking heuristic
when the required time becomes comparable to (or greater than) the
required time for the full, untweaked process. In our experiments, one
tweaking iteration is 5-10 times faster than one \SLP iteration and
one Newton step on \eqref{eq:KKT.LS}.  For each LP subproblem, the
dual simplex option is used in \CPLEX.  For \IPOPT and MIPS, all
parameters are set to their default values.  In the experiments, the
voltage limits on each system are chosen so that the optimal voltage
magnitude values for the unmodified problems do not violate the
voltage constraints.

\subsection{IEEE 57-Bus System with Line Impedance Perturbations}
\label{sec:57}

\begin{table}
\centering
\subfloat[IEEE 57-Bus System for Line Impedance Perturbation.]{
\begin{tabular}{||c|c|c|c|c|c|c|c|c|c|c||}\hline
  \mc{3}{||c|}{\# of Buses}
  & \# of
  & \mc{4}{c|}{\# of Variables}
  & \mc{2}{c|}{\# of Constraints}
  & \mr{2}{*}{$(\lb{V},\ \ub{V})$} \\
\cline{1-3}\cline{5-10}
  Ref & $\cG$ & $\cD$
  & Lines
  & $V$ & $\theta$ & $\sigma^{\pm},\rho$ & Total
  & $=$ & $\le$ or $\ge$
  & \\\hline\hline
1 & 6 & 50 & 80 & 56 & 50 & 37 & 143 & 106 & 186 & (0.93,\ 1.07) \\\hline
\end{tabular}\label{tbl:IEEE57}}

\subfloat[\SLP Results on IEEE 57-Bus System.]{
\begin{tabular}{||c||c|c|c|c|c|c|c||}\hline
  \mr{2}{*}{$\beta$}
  & \mr{2}{*}{$n-m-\abs{\cA^*}$}
  & \mc{2}{c|}{Iterations} & Time & \mc{2}{c|}{Loads Shed} & $\#$ of Buses w/\\
\cline{3-4}\cline{6-7}
    &   & \SLP & AS & (sec) & $P$ & $Q$ & Loads Shed \\\hline\hline
1.0 & 0 & 3(0) & -  & 0.03  &     0  &    0 & 0 \\\hline
1.2 & 0 & 3(2) & -  & 0.04  &  2.93  & 1.46 & 2 \\\hline
1.4 & 0 & 3(5) & -  & 0.04  &  8.37  & 3.80 & 4 \\\hline
1.6 & 2 & \tcell{23(35)\\5(20)} & \tcell{-\\3(0)} & \tcell{0.13\\0.06} & 18.06 &  7.89 &  9 \\\hline
1.8 & 2 & \tcell{23(47)\\7(20)} & \tcell{-\\3(1)} & \tcell{0.14\\0.8}  & 27.14 & 12.57 & 10 \\\hline
2.0 & 2 & \tcell{28(72)\\6(31)} & \tcell{-\\3(2)} & \tcell{0.15\\0.05} & 35.65 & 16.57 & 11\\\hline
\mc{8}{r}{($P$ in MW / $Q$ in MVAr)}
\end{tabular}\label{tbl:IEEE57.result}}
\caption{Line Impedance Perturbation: IEEE 57-Bus System.}
\end{table}


We study a standard test case, the IEEE 57-Bus system, to explore the
overall behavior of our algorithm. Information about the system is
provided in Table~\subref*{tbl:IEEE57}. For the purpose of obtaining
power network instances which are disrupted in different degrees of
severity from the 57-Bus systems data, we modify the specifications of
the problem by multiplying the impedance of each line by factor
$\beta$ in the range $[1,2]$. (Although disruptions of this kind would
not happen in practice, we believe that this technique is a reasonable
way to define a sequence of increasingly stressed grids that are
closely related to the realistic grids that are found in standard test
sets.) The values given in \MATPOWER data file are used as the initial
starting point for the algorithm.

Table~\subref*{tbl:IEEE57.result} shows the results for the each value
of $\beta$.  The column ``$n-m-\abs{\cA^*}$'' indicates whether the
solution is fully determined by its constraints (the value is zero) or
not (a positive integer).  For the fully-determined cases, the active
set heuristic is often not invoked, in which case we show the results
as a single line. In cases for which the active-set heuristics is
used, we show two lines of results, one for the case in which this
heuristic is turned off and one for the case in which we allow it to
be invoked. The column labelled \SLP shows the number of \SLP
iterations with the total number of dual simplex iterations in
parentheses and the column labelled AS shows the number of active-set
iterations with the total number of tweaking iterations in
parentheses. (We use the same format to present the number of
iterations for all the following tables.)

When $\beta=1.0$ (the original, non-disrupted system), the power flow
problem is feasible and thus the algorithm takes the same steps as
Newton's method applied to the power-flow equations. It converges {\em
without} needing any simplex iterations; the original basis
factorization is enough.  As expected, the solution does not require
any load shedding, and the objective is therefore zero.  With
$\beta=1.2$ or $1.4$, the solution is still fully determined by the
active constraints, and we observe quadratic convergence, consistently
with the analysis of Subsection~\ref{sec:slp.local}, though in these
cases, load shedding is required. A few simplex iterations are needed
in the process of solving the linear programming subproblems, to
resolve the active sets.

The cases of $\beta=1.6$, $1.8$, and $2.0$ show the benefits of the
active set heuristic.
When $\beta=2.0$, for example, the \SLP algorithm without active-set
heuristic requires 28 iterations to find the solution, converging to
this underdetermined solution at a slow linear rate.
When the active-set heuristic is used, the algorithm invokes it after
six iterations of \SLP and converges rapidly thereafter.  The initial
active set estimation from the \SLP steps was not exact for this case,
but two iterations of the active-set heuristic sufficed to find the
optimal active set.


\begin{table}\centering
\begin{tabular}{||c||c||c|c||c|c||c|c|c||c||}
\hline
\mr{3}{*}{$\beta$}
  & \mr{3}{*}{\tcellr{Bus\\No.}}
  & \mc{2}{c||}{Demanded}
  & \mc{2}{c||}{Injected}
  & \mc{3}{c||}{Loads Shed}
  & Buses \\
\cline{3-9}
  &
  & \mr{2}{*}{$P$} & \mr{2}{*}{$Q$} & \mr{2}{*}{$P$} & \mr{2}{*}{$Q$}
  & \mr{2}{*}{\%} & \mc{2}{c||}{Total}
  & with \\
\cline{8-9}
  &
  & & & &
  & & $P$ & $Q$
  & $V=\lb{V}$ \\
\hline
\hline
 1.6 & \tcellr{20\\30\\31\\32\\33\\42\\53\\56\\57}
     & \tcellr{2.3\\3.6\\5.8\\1.6\\3.8\\7.1\\20\\7.6\\6.7}
     & \tcellr{1.0\\1.8\\2.9\\0.8\\1.9\\4.4\\10\\2.2\\2.0}
     & \tcellr{2.02\\3.52\\0\\1.26\\0\\5.39\\19.00\\6.98\\3.26}
     & \tcellr{0.88\\1.76\\0\\0.63\\0\\3.34\\9.50\\2.02\\0.97}
     & \tcellr{12.3\%\\2.0\%\\100.0\%\\21.2\%\\100\%\\24.0\%\\5.0\%\\8.2\%\\51.3\%}
     & 17.06
     & 7.89
     & \tcellr{20\\26\\34\\42\\53\\56\\57} \\
\hline
 2.0 & \tcellr{20\\25\\30\\31\\32\\33\\35\\42\\53\\56\\57}
     & \tcellr{2.3\\6.3\\3.6\\5.8\\1.6\\3.8\\6.0\\7.1\\20\\7.6\\6.7}
     & \tcellr{1.0\\3.2\\1.8\\2.9\\0.8\\1.9\\3.0\\4.4\\10\\2.2\\2.0}
     & \tcellr{0.11\\5.29\\0.09\\0\\0.30\\0\\4.27\\3.71\\14.53\\4.78\\2.08}
     & \tcellr{0.05\\2.69\\0.05\\0\\0.15\\0\\2.14\\2.30\\7.26\\1.38\\0.62}
     & \tcellr{95.3\%\\16.1\%\\97.5\%\\100.0\%\\81.2\%\\100.0\%\\28.8\%\\47.7\%\\27.4\%\\37.1\%\\69.0\%}
     & 35.65
     & 16.57
     & \tcellr{20\\26\\35\\42\\53\\56\\57} \\
\hline
\mc{10}{r}{($P$ in MW / $Q$ in MVAr)}
\end{tabular}
\caption{Feasibility Restoration Results on IEEE 57-Bus Systems with
  $\beta=1.6$ and $\beta=2.0$}\label{tbl:fr_ieee57}
\end{table}

When applied to these problems with $\beta$ between $1.0$ and $1.8$,
the \MATPOWER AC power flow solver obtained solutions using Newton's
method within ten iterations. Thus, there exist solutions to the basic
power flow equations \eqref{eq:intro.pflow} without load shedding for
these cases --- but the voltage magnitudes at some buses are too low
for these solutions to be practically operational. When $\beta=2.0$,
\MATPOWER fails to find a solution, suggesting strongly that the
power-flow equations do not have a solution, even one with impractical
voltage magnitudes.

Details of the solutions obtained by the \SLPAS algorithm for the
cases $\beta=1.6$ and $\beta=2.0$ are shown in
Table~\ref{tbl:fr_ieee57}.  Load shedding is required at 9 and 11
buses, respectively, the table showing the percentage of load shedding
required at each bus. The last column indicates those buses at which
the optimal voltage magnitudes are at their lower bounds.

\subsection{Loss of a Transmission Line (N-1 Cases)}

A more practical way to disrupt a power system is to remove
transmission lines. We test our algorithm on the standard IEEE 300-Bus
system and Polish systems by removing each transmission line in turn
--- the so-called ``$N-1$'' cases. Since we are mainly interested in
unsolvable cases, we only consider the lines such that ({\em a}) the
whole system remains fully connected when the line is removed; and
({\em b}) the Newton's method of \MATPOWER fails to solve the
power-flow equations for the damaged system.  Using these criteria, we
found sixteen, two, and two-hundred sixty eight unsolvable cases for
IEEE 30-Bus system, Polish 2383-Bus system, and Polish 2746-Bus
system, respectively. However, not all unsolvable cases have a
load-shedding solution that can make the disrupted system feasible.
We applied our algorithm on the unsolvable instances and show the
result of the cases on which there exists a load-shedding solution in
Table~\ref{tbl:n1}.

\begin{table}
\centering
\subfloat[IEEE 300-Bus, Polish 2383-Bus and Polish 2746-Bus Systems for $N-1$ Cases.]{
\begin{tabular}{||c||c|c|c|c|c|c||}
\hline
\mr{2}{*}{System}
  & \# of
  & \# of
  & \# of
  & \mc{2}{c|}{\# of Constraints}
  & \mr{2}{*}{$(\lb{V},\ \ub{V})$} \\ 
\cline{5-6}
  & Buses
  & Lines
  & Variables
  & \makebox[25pt]{$=$} & $\le$ or $\ge$ 
  & \\ 
\hline
\hline
IEEE 300    &  300 &  410 &  757 &  530 &  916 & (0.92,\ 1.08) \\\hline
Polish 2383 & 2383 & 2895 & 6087 & 4438 & 7410 & (0.90,\ 1.12) \\\hline
Polish 2746 & 2746 & 3278 & 6925 & 5127 & 8360 & (0.98,\ 1.20) \\\hline
\end{tabular}\label{tbl:n1.info}}


\subfloat[\SLPAS Results on IEEE 300-Bus System]
{\begin{tabular}{||c||c|c|c|c|c|c||}
\hline
Line
  & \mc{2}{c|}{Iterations}
  & Time & \mc{2}{c|}{Loads Shed} & \# of Buses w/\\
\cline{2-3}\cline{5-6}
Removed
  & \SLP & AS
  & (sec) & $P$ (MW) & $Q$ (MVAr) & Loads Shed \\
\hline
\hline
 66 & 8(9) & 3(0) & 0.08 & 49.88 & 13.80 & 2\\ \hline
114 & 5(1) & - & 0.05 &  136.63 & 66.43 & 1 \\ \hline
177 & 5(5) & - & 0.05 &    8.15 & 375.82 & 1 \\ \hline
181 & 6(7) & - & 0.05 & 1089.30 & 153.88 & 3\\ \hline
182 & 5(6) & - & 0.05 &  246.74 & 230.11 & 2 \\ \hline
268 & 3(10) & 4(0) & 0.07 & 1158.70 & 39.06  & 4\\ \hline
294 & 5(10) & - & 0.05 & 154.37 & 0.03 & 3\\ \hline
309 & 9(30) & 6(4) & 0.13 & 822.25 & 21.23 & 7\\ \hline
364 & 5(1) & - & 0.05 & 107.78 &  50.83 & 1\\ \hline
367 & 5(4) & - & 0.05 & 229.10 & 114.12 & 2\\ \hline
369 & 6(8) & - & 0.05 &  51.04 &   4.79 & 2\\ \hline
370 & 5(1) & - & 0.05 &  34.81 &   2.49 & 1\\ \hline
381 & 5(4) & - & 0.06 &  23.27 &   7.96 & 1\\ \hline\hline
\mc{7}{||c||}{\SLPAS did not converge when line 116, 187 or 350 is removed}\\\hline
\end{tabular}\label{tbl:n1.ieee300}}

\subfloat[\SLPAS Results on Polish 2383-Bus System]
{\begin{tabular}{||c||c|c|c|c|c|c||}
\hline
Line
  & \mc{2}{c|}{Iterations}
  & Time & \mc{2}{c|}{Loads Shed} & \# of Buses w/\\
\cline{2-3}\cline{5-6}
Removed
  & \SLP & AS
  & (sec) & $P$ (MW) & $Q$ (MVAr) & Loads Shed \\
\hline
\hline
466 & 6(29) & - & 0.15 & 183.70 & 33.61 & 17\\ \hline
469 & 8(35) & 3(0) & 0.38 & 183.78 & 33.36  & 18\\ \hline
\end{tabular}\label{tbl:n1.polish2383}}

\subfloat[\SLPAS Results on Polish 2746-Bus System: 8 out of 160 Solved Instances Shown]
{\begin{tabular}{||c||c|c|c|c|c|c||}
\hline
Line
  & \mc{2}{c|}{Iterations}
  & Time & \mc{2}{c|}{Loads Shed} & \# of Buses w/\\
\cline{2-3}\cline{5-6}
Removed
  & \SLP & AS
  & (sec) & $P$ (MW) & $Q$ (MVAr) & Loads Shed \\
\hline
\hline
  28 & 4(48)& - & 0.17 & 308.77 &  58.25  & 1 \\ \hline
 104 & 4(3) & - & 0.10 & 223.04 & 76.02 & 1 \\ \hline
1458 & 4(3) & - & 0.10 &   0.00 &  0.00 & 0 \\ \hline
1957 & 4(8) & - & 0.11 &   4.28 &  0.00 & 1 \\ \hline
2068 & 4(55) & - & 0.23 &  25.30 &   4.69  & 1 \\ \hline
3037 & 4(8) & - & 0.11 &  12.47 &  1.50 & 3 \\ \hline
3046 & 4(8) & - & 0.11 &   9.62 &  0.76 & 2 \\ \hline
3496 & 4(51) & - & 0.18 &   1.53 &   0.50  & 1 \\ \hline
\mc{7}{||c||}{\SLPAS solved 160 instances, and did not converge on 130 infeasible instances.}\\\hline
\end{tabular}\label{tbl:n1.polish2746}}
\caption{Loss of a Transmission Line: IEEE 300-Bus System and Polish Systems.}\label{tbl:n1}
\end{table}

Details of the $N-1$ systems derived from the IEEE 300-Bus system and
the Polish systems are presented in Table~\subref*{tbl:n1.info}. (The
number of lines listed in the table is less than the number of lines
in the original system by one. For example, the IEEE 300-Bus system
has 411 transmission lines, but since one of them is removed in each
of our cases, we list 410 lines.)

Table~\subref*{tbl:n1.ieee300} shows results for the sixteen
unsolvable $N-1$ cases. As mentioned, none of these cases can be
solved by the Newton's method, and three of them cannot even be solved
as feasibility restoration problems, probably because of the
restrictions imposed in our formulation \eqref{eq:feas.rest}. We
verified that these three cases are actually infeasible problems
using a different nonlinear programming solver. For the ten cases in
which the solutions are fully determined by the constraints, we see
fast convergence of \SLP without the need to invoke the active-set
heuristic. In the three underdetermined cases, the active-set
heuristic is used for speedup.

We also solved the feasibility restoration problems formulated for the
two $N-1$ unsolvable cases of Polish 2383-Bus successfully using the
\SLPAS algorithm; Table~\subref*{tbl:n1.polish2383} shows the results.
For Polish 2746-Bus system, the \SLPAS algorithm solves 160 out of 290
instances, and did not converge for 130 instances, which are confirmed
to be infeasible using a different solver.  In
Table~\subref*{tbl:n1.polish2746}, we present the results of eight
typical instances, which are solved by our algorithm. We note that no
load shedding is required when transmission line 1458 is removed,
which means the $N-1$ instance is feasible. (\MATPOWER's Newton solver
fails on this instance because the Jacobian of the power flow
constraints goes singular.)

Note that in all cases, feasibility can be restored to the system by
shedding load on just a few buses if needed.

\subsection{Loss of a Generator}
Another possible contingency that can disrupt a power system is a
generator failure. We simulate this situation by disabling a
generator. Similarly to the previous subsection, we identified seven
unsolvable cases by removing one generator at a time from IEEE 300-Bus
system for all generators (thus, the corresponding bus is converted
into a demand bus in the AC power flow problem), and formulated the
feasibility restoration problem for these cases. The problem setup
(see Table~\subref*{tbl:n1.gen.info}) is also similar to the previous
experiments (cf. Table~\subref*{tbl:n1.info}), but the number of
variables and the number of lines are changed since a generator is
removed (and the node to which the generator attached is converted
into a demand bus) instead of a transmission line. Since the Polish
2383-Bus system and Polish 2746-Bus system do not have any unsolvable
instances caused by disabling a generator, they are not tested in this
experiment.

Table~\subref*{tbl:n1.gen.ieee300} shows that the \SLPAS algorithm
successfully solves all seven cases which have a loss of a generator,
and the active-set heuristic is invoked for six instances.

\begin{table}
\centering
\subfloat[IEEE 300-Bus System for Loss of A Generator.]{
\begin{tabular}{||c||c|c|c|c|c|c||}
\hline
\mr{2}{*}{System}
  & \# of
  & \# of
  & \# of
  & \mc{2}{c|}{\# of Constraints}
  & \mr{2}{*}{$(\lb{V},\ \ub{V})$} \\ 
\cline{5-6}
  & Buses
  & Lines
  & Variables
  & \makebox[25pt]{$=$} & $\le$ or $\ge$ 
  & \\ 
\hline
\hline
IEEE 300    &  300 &  411 &  756 &  530 &  916 & (0.92,\ 1.08) \\\hline
\end{tabular}\label{tbl:n1.gen.info}}

\subfloat[\SLPAS Results on IEEE 300-Bus System]
{\begin{tabular}{||c||c|c|c|c|c|c||}
\hline
Generator
  & \mc{2}{c|}{Iterations}
  & Time & \mc{2}{c|}{Loads Shed} & \# of Buses w/\\
\cline{2-3}\cline{5-6}
Disabled
  & \SLP & AS
  & (sec) & $P$ (MW) & $Q$ (MVAr) & Loads Shed \\
\hline
\hline
 98 &  16(89) & 5(10) & 0.19 & 1240.16 & 603.89 & 14 \\\hline
165 &   6(73) &  4(3) & 0.10 &  409.36 &  13.12 & 19 \\\hline
166 &   6(69) &  4(3) & 0.10 &  411.60 &  13.12 & 19 \\\hline
170 &  9(219) &  1(0) & 0.09 & 1377.83 & 142.14 & 15 \\\hline
249 &   5(76) &  4(3) & 0.10 &  353.79 &  11.47 & 19 \\\hline
252 &   5(1)  &    -  & 0.05 &  375.55 & 147.28 &  1 \\\hline
264 &   5(82) &  4(1) & 0.08 &  436.88 &  11.87 & 21\\ \hline
\mc{7}{r}{Generator number is as appeared in the \MATPOWER data file (\texttt{case300.m})}
\end{tabular}\label{tbl:n1.gen.ieee300}}
\caption{Loss of A Generator: IEEE 300-Bus System.}
\end{table}

\subsection{Comparisons with \IPOPT and MIPS}\label{sec:performance.comp}
\begin{sidewaystable}
\centering
\subfloat[IEEE Standard Test Cases and Polish Systems for Performance Comparison.]{
\makebox[.8\textwidth]{
\begin{tabular}{||c||c|c|c|c|c|c||}
\hline
\mr{2}{*}{System}
  & \# of
  & \# of
  & \# of
  & \mc{2}{c|}{\# of Constraints}
  & \mr{2}{*}{$(\lb{V},\ \ub{V})$}\\
\cline{5-6}
  & Buses 
  & Lines
  & Variables
  & $=$ & $\le$ or $\ge$ & \\\hline\hline
IEEE 118 & 118 & 186 & 262 & 181 & 290 & (0.93, 1.07) \\\hline
IEEE 300 & 300 & 411 & 757 & 530 & 916 & (0.92, 1.08) \\\hline
Polish 2383 & 2383 & 2896 & 6087 & 4438 & 7410 & (0.90, 1.12) \\\hline
Polish 2746 & 2746 & 3279 & 6925 & 5127 & 8360 & (0.98, 1.20) \\\hline
\end{tabular}\label{tbl:fr_performance}}}

\subfloat[\SLPAS, \IPOPT and MIPS Results on IEEE Standard Test Cases and Polish Systems.]{
\begin{tabular}{||c|c||c|c|c||c|c||c|c||c|c||c|c|c||}
\hline
\mr{3}{*}{System}
  & \mr{3}{*}{$\beta$}
  & \mc{3}{c||}{\SLPAS}
  & \mc{4}{c||}{\IPOPT}
  & \mc{2}{c||}{MIPS}
  & \mc{2}{c|}{Loads Shed}
  & \# of Buses\\
\cline{3-5}\cline{6-11}\cline{12-13}
  &
  & \mc{2}{c|}{Iterations}
  & \mr{2}{*}{Time (s)}
  & \mc{2}{c||}{L-BFGS}
  & \mc{2}{c||}{Hessian}
  & \mr{2}{*}{Iter.}
  & \mr{2}{*}{Time (s)}
  & $P$ & $Q$ & with \\
\cline{3-4}\cline{6-9}
  &
  & LP & AS &
  & Iter. & Time (s) & Iter. & Time (s)
  & &
  & (MW) & (MVAr) & Loads Shed \\
\hline
\hline
\mr{4}{*}{\tcell{118 Bus}}
  & 1.5 & 4(0)  & -    & 0.03 &  6 & 0.03 &  7 & 0.03 & 12 & 0.05 &      0 &     0 & 0\\
  & 2.0 & 4(2)  & -    & 0.03 &  9 & 0.05 & 10 & 0.04 & 13 & 0.05 &  10.54 &  5.53 & 2\\
  & 2.5 & 5(15) & 3(1) & 0.06 & 15 & 0.08 & 16 & 0.06 & 16 & 0.06 &  62.81 & 25.67 & 9\\
  & 3.0 & 7(20) & 3(0) & 0.06 & 14 & 0.08 & 15 & 0.06 & 16 & 0.06 & 178.21 & 70.18 & 15\\
\hline
\mr{2}{*}{\tcell{300 Bus}}
  & 1.1 &  5(5) & - & 0.05 & 13 & 0.12 & 18 & 0.14 & 19 & 0.16 &  38.39 &  23.54 & 5\\
  & 1.2 & 5(11) & - & 0.04 & 18 & 0.18 & 18 & 0.13 & 19 & 0.16 & 222.54 & 100.13 & 11\\
\hline
\mr{5}{*}{\tcell{2383 Bus}}
  & 1.2 & 8(27)   &  -   & 0.24 &  24 & 1.29 & 20 & 0.91 & 22 & 1.76 & 188.66 &  41.44 & 17\\
  & 1.4 & 11(51)  & 5(0) & 0.65 &  23 & 1.25 & 20 & 0.91 & 21 & 1.64 & 233.97 &  55.02 & 22\\
  & 1.6 & 6(82)   & 5(4) & 0.62 &  26 & 1.47 & 21 & 0.95 & 21 & 1.65 & 290.20 &  67.03 & 29\\
  & 1.8 & 14(283) & 5(8) & 0.96 &  46 & 2.61 & 31 & 1.38 & 38 & 3.05 & 535.43 &  87.65 & 52\\
  & 2.0 & 13(571) & 5(2) & 0.89 &  95 & 5.99 & 44 & 1.46 & 44 & 2.65 & 855.74 & 116.32 & 82\\
\hline
\mr{5}{*}{\tcell{2746 Bus}}
  & 1.2 &   5(1) & 1(0) & 0.22 & 12 & 0.73 & 12 & 0.75 & 16 & 1.31 &   3.52 &   3.08 & 1\\
  & 1.4 &   5(5) & -    & 0.14 & 21 & 1.36 & 20 & 1.12 & 19 & 1.57 &  42.87 &  16.93 & 7\\
  & 1.6 &  6(23) & 4(2) & 0.52 & 38 & 2.59 & 28 & 1.54 & 26 & 2.19 & 176.44 &  51.98 & 20\\
  & 1.8 &  9(61) & 3(0) & 0.54 & 40 & 2.76 & 26 & 1.42 & 27 & 2.24 & 337.52 &  83.29 & 35\\
  & 2.0 & 8(182) & 4(4) & 0.70 & 55 & 3.75 & 38 & 2.07 & 31 & 2.59 & 590.98 & 112.53 & 55\\
\hline
\end{tabular}}
\caption{Performance Comparisons: \SLP (with Active-Set Heuristic) vs
  \IPOPT and MIPS.}\label{tbl:fr_perform}
\end{sidewaystable}

Performance comparisons between \SLPAS (that is, \SLP with the
active-set heuristic) and the solvers \IPOPT and MIPS for feasibility
restoration problems are shown in Table~\ref{tbl:fr_perform}. We used
four test problems from \MATPOWER{}, which include the IEEE standard
test cases and Polish systems, with details shown in
Table~\subref*{tbl:fr_performance}.  We used a similar setup here to
Subsection~\ref{sec:57}, modifying the basic IEEE data and Polish
system data by scaling the impedance of all lines by a factor $\beta >
1$.

Two different variants of \IPOPT are tested. One setting uses the
exact Hessian of the Lagrangian, while the other uses a limited-memory
quasi-Newton (L-BFGS) approximation.
MIPS, which is provided within \MATPOWER, uses the exact Hessian, as
does the active-set heuristic within \SLPAS.

The 118-Bus system is quite robust. Even with $\beta=3$, there exists
a solution that can operate the system within the given voltage range
with load adjustments on just a few buses. All methods solve the
problem in less than .1 seconds.

On the 300-Bus system, the solutions were fully determined in the
cases reported. The \SLPAS converges rapidly without needing to invoke
the active-set heuristic, and gives faster runtimes than rival
approaches. We also tried setting $\beta=1.3$, but none of the
algorithms could find a solution for this case.

For the 2383-Bus systems, \SLPAS is slightly faster than the other
approaches, but its advantage over \IPOPT becomes slimmer as $\beta$
is increased.  On these heavily disrupted systems, many buses require
adjustments to recover feasibility, and the solutions are
underdetermined, so the dual simplex algorithm requires many more
iterations on each subproblem. The active-set heuristic continues to
work well, however, and yields fast local convergence.

\SLPAS is significantly faster for the cases based on the 2476-bus
system, even though there are underdetermined solutions in some of
these cases.

We conclude with some observations about the implementations of each
of these solvers. Note that both \IPOPT (with exact Hessian) and MIPS
both use interior point approaches and require similar numbers of
iterations, but \IPOPT is about 30\% faster on the more difficult
problems, probably because it is implemented in \MATLAB whereas \IPOPT
is implemented in \verb!C++!.  Our \SLPAS algorithm is coded in a
combination of \MATLAB (to set up the subproblems and to perform the
active-set heuristic) and native code (the \CPLEX solver for the LP
subproblems) and thus falls somewhere between \IPOPT and MIPS in the
efficiency of its code. A fully native implementation of \SLPAS has
the potential to perform even better.

\section{Conclusions} \label{sec:con}

We have proposed an algorithm consisting of \SLP and an active-set
heuristic to solve nonsmooth penalty-function formulations of
nonlinear programming problems. Global and local convergence
properties were explored, with a focus on local (quadratic)
convergence in the case in which the solution is fully determined by
the constraints --- a situation that occurs frequently in our target
application. In the final sections, we described application of our
algorithm to the problem of restoring feasibility to a disrupted power
system, in a practical way that limits the number of buses at which
load must be shed. Computational results were presented for system
sizes up to 2746 buses, and comparisons were performed with
interior-point solvers on the same formulations.

Since the \SLP algorithm needs to solve LP subproblems at each
iteration, it may become slow due to the overhead of finding the LP
solutions for large problems. For highly perturbed power systems that
require many demand buses to be adjusted to recover feasibility, many
simplex iterations may be required to find a subproblem solution.
Because of the nonconvex formulation, our approach is not guaranteed
to find a global solution, a limitation shared with other solvers,
including interior-point methods.

One goal of our formulation and method is to guide system operators
toward load-shedding patterns that restore practical operation of the
grid at minimum disruption. Another goal is to use the optimal
objective in \eqref{eq:feas.rest} as a measure of disruption to the
grid, to be used in analyzing the vulnerability of the grid to
deliberate attacks or natural disturbances. Vulnerability analysis may
indicate what capital improvements could make the system
more robust to such disruptions. We are exploring these issues further
in current research.


\appendix

\section{Solving Augmented Linear Systems}\label{appdx:AugLS}

Assume that we already know the solution $x_0$ of the square linear
system
\[ Hx = b_1, \]
where $H$ is symmetric, as well as factors $L$ (lower triangular) and
$D$ (block diagonal) such that $LDL^T=H$. Suppose we are presented
with the following augmented square linear system:
\[ \bm{H & V\\X & S}\bm{x_1\\ x_2} = \bm{b_1\\b_2}, \]
where the dimensions of the square matrix $S$ are much smaller than
those of $H$.  We first factorize the matrix into the two block
triangular matrices:
\[ \bm{H & V\\ X & S} = \bm{ H & 0 \\ X & C}\bm{ I & Y \\ 0 & I } \]
where the matrices $C$ and $Y$ can be obtained by solving the following
problems:
\begin{alignat*}{2}
  V &= HY   & \quad\Rightarrow\quad & LDL^TY = V, \\
  S &= XY+C & \quad\Rightarrow\quad & C = S-XY.
\end{alignat*}
Because of the properties of $L$ and $D$, the matrix $Y$ can be
calculated economically, while $C$ requires simply a matrix
multiplication.  We can therefore rewrite the augmented system as two
linear systems with auxiliary variables $w_1$ and $w_2$:
\begin{alignat*}{2}
  \bm{H & 0 \\ X & C} \bm{w_1\\ w_2}         &= \bm{b_1\\ b_2}, &
  \qquad  \bm{I & Y\\ 0  & I} \bm{x_1\\ x_2} &= \bm{w_1\\ w_2}.
\end{alignat*}
Since $Hw_1=b_1$, we have $w_1=x_0$. Then $w_2$ is obtained by solving
the following system:
\[ Cw_2 = b_2-Xw_1, \]
which can be performed economically, since $C$ is small.  We can then
find $x_1$ and $x_2$ by setting
\begin{align*}
 x_2 &= w_2, \\
 x_1 &= w_1-Yx_2.
\end{align*}

\bibliographystyle{spbasic}
\bibliography{feasibility}

\end{document}